\documentclass[a4paper,12pt]{article}

\usepackage[utf8]{inputenc}
\usepackage[mathscr]{eucal}
\usepackage{amsmath}
\usepackage{amssymb}
\usepackage{MnSymbol}
\usepackage{amsfonts}
\usepackage{graphicx}
\usepackage{multicol}
\usepackage{enumerate}

\newcommand{\df}[1] {{\bfseries  #1}}

\newcommand{\zbp}{\varnothing}
\newcommand{\R}{\mathbb{R}}

\newcommand{\Z}{\mathbb{Z}}

\newcommand{\vare}{\varepsilon}

\newcommand{\dom}{\text{dom}}

\newcommand{\add}{\text{add}}
\newcommand{\cov}{\text{cov}}
\newcommand{\cof}{\text{cof}}
\newcommand{\non}{\text{non}}

\newcommand{\diam}{\text{diam}}

\newcommand{\oA}{{\mathcal{A}}}
\newcommand{\oB}{{\mathcal{B}}}

\newcommand{\oI}{{\mathcal{I}}}

\newcommand{\oM}{{\mathcal{M}}}
\newcommand{\oN}{{\mathcal{N}}}
\newcommand{\oO}{{\mathcal{O}}}
\newcommand{\oU}{{\mathcal{U}}}
\newcommand{\oV}{{\mathcal{V}}}
\newcommand{\oW}{{\mathcal{W}}}
\newcommand{\oF}{{\mathcal{F}}}

\newcommand{\oS}{{\mathcal{S}}}
\newcommand{\oX}{{\mathcal{X}}}

\newcommand{\bont}{{\mathfrak{b}}}

\newcommand{\dont}{{\mathfrak{d}}}
\newcommand{\Kont}{{\mathfrak{K}}}

\newcommand{\QQ}{\boldsymbol{Q}}

\newcommand{\SN}{\boldsymbol{\text{S}\oN}}
\newcommand{\MM}{\boldsymbol{\oM}}
\newcommand{\PM}{\boldsymbol{\text{P}\oM}}
\newcommand{\PkM}{\boldsymbol{\text{P}_\kappa\oM}}

\newcommand{\NS}{\boldsymbol{\oN\oS}}
\newcommand{\bzero}{\boldsymbol{0}}
\newcommand{\bjeden}{\boldsymbol{1}}

\newcommand{\obc}{\mathord{\upharpoonright}}
\newcommand{\Spl}{\text{Split}}

\newcommand{\Lim}{\text{Lim}}

\newcommand{\len}{\text{len}}

\newtheorem{theo}{Theorem}
\newtheorem{lemma}[theo]{Lemma}
\newtheorem{prop}[theo]{Proposition}
\newtheorem{cor}[theo]{Corollary}

\newtheorem{prob}[theo]{Problem}

\title{Special subsets of the generalized Cantor space $2^\kappa$ and generalized Baire space $\kappa^\kappa$}
\author{Michał Korch, Tomasz Weiss}
\date{}

\begin{document}

\maketitle

\begin{abstract}
  In this paper, we are interested in parallels to the classical notions of special subsets in $\R$ defined in the generalized Cantor and Baire spaces ($2^\kappa$ and $\kappa^\kappa$). We consider generalizations of the well-known classes of special subsets, like Lusin sets, strongly null sets, concentrated sets, perfectly meagre sets, $\sigma$-sets, $\gamma$-sets, sets with Menger, Rothberger or Hurewicz property, but also of some less-know classes like $X$-small sets, meagre additive sets, Ramsey null sets, Marczewski, Silver, Miller and Laver-null sets. We notice that many classical theorems regarding these classes can be relatively easy generalized to higher cardinals although sometimes with some additional assumptions. This paper serves as a catalogue of such results along with some other generalizations and open problems.
\end{abstract}

\section{Introduction and preliminaries}

Many classical notions of special subsets of $2^\omega$ can be transferred to the case of the generalized Cantor space $2^\kappa$. In this paper we study those classes of sets in this setting. In some cases, when it seems more appropriate, we also study such classes in the generalized Baire space $\kappa^\kappa$. It turns out that many properties of subsets of $2^\omega$ or of $\omega^\omega$ can be easily proved in $2^\kappa$ or $\kappa^\kappa$, although sometimes one has to use some additional set-theoretic assumptions. Next we deal with less common classes of small sets in $2^\kappa$ and $\kappa^\kappa$.

In most cases the theorems presented in this paper are simple generalizations of the classical theorems, although sometimes we need some additional assumptions. In such cases we make an attempt to discuss these assumptions or state some open problems. In particular, such considerations can be found in sections~\ref{s-sn}, \ref{s-pn}, \ref{s-rn}, and \ref{s-m0}.

\subsection{Trees}

Fix any set $A$ and an ordinal number $\xi$. Given a sequence $t\in A^\alpha$ with $\alpha<\xi$, we denote $\alpha=\len(t)$. A~set $T\subseteq A^{<\xi}$ will be called a~\df{tree} if for all $t\in T$ and $\alpha< \len(t)$, $t\obc \alpha\in T$ as well, and for every $\alpha<\xi$, there exists $t\in T$ such that $\len(t)=\alpha$. A~branch in a~tree is a~maximal chain in it. For a~tree $T\subseteq A^{<\xi}$ of height $\xi$, let $[T]=\{x\in A^\xi\colon \forall_{\alpha<\xi} x\obc\alpha\in T\}$.

A~node $s\in T\subseteq A^{<\xi}$ is called a~\df{branching point} of $T$ if $s^\frown a,s^\frown b\in T$ for some distinct $a,b\in A$. The set of all branching points of a~tree $T$ is denoted by $\Spl(T)$. For $\alpha<\xi$, $t\in \Spl_\alpha(T)$ if $\langle\{s\subsetneq t\colon s\in \Spl(T)\},\subseteq\rangle$ is order isomorphic with $\alpha$. 

A~tree $T\subseteq A^{<\xi}$ is \df{perfect} if for any $t\in T$, there exists $s\in T$ such that $t\subseteq s$ and $s\in\Spl(T)$. 

A~tree $T\subseteq A^{<\xi}$ is pruned, if its every maximal chain is of length $\xi$. 

Notice that a set $C\subseteq A^\omega$ is closed if an only if $C=[T]$ for a pruned tree $T$. We denote such tree by $T_C$. Moreover, a~set $P\subseteq 2^\omega$ is perfect if and only if $T_P$ is a~perfect tree. Notice also that a closed set $C\subseteq \omega^\omega$ is compact if and only if there exists a sequence $\langle n_i\rangle_{i\in\omega}$ such that if $x\in C$, then $x(i)<n_i$ for all $i\in\omega$. 
 
A perfect tree $T\subseteq A^{<\omega}$ is called a \df{Silver perfect tree} if \[\forall_{w,v\in T}\left( \len(v)=\len(w)\Rightarrow\forall_{j\in A}(w^\frown j\in T\Rightarrow  v^\frown j\in T)\right).\]

A perfect tree $T\subseteq \omega^{<\omega}$ is called a~\df{Laver perfect tree} if there exists $s\in T$ such that for all $t\in T$, either $t\subseteq s$, or $\left|\left\{n\in\omega\colon t^\frown  n\in T\right\}\right|=\aleph_0$.

Similarly, a perfect tree $T\subseteq \omega^{<\omega}$ is called a~\df{Miller perfect tree} if for every $s\in T$ there exists $t\in T$ such that  $s\subseteq t$, and $\left|\left\{n\in\omega\colon t^\frown n\in T\right\}\right|=\aleph_0$.

A set $P\subseteq 2^{<\omega}$ is called \df{Silver perfect set} if  $T_P$ is a Silver perfect tree. Analogously, a set $P\subseteq \omega^\omega$ is called \df{Laver (respectively, Miller) perfect set} if $T_P$ is a Laver (respectively, Miller) perfect tree.

\subsection{Introducing the generalized Cantor space $2^\kappa$ and the generalized Baire space $\kappa^\kappa$}\label{intro-gen}
In this paper, we consider the generalized Cantor space $2^\kappa$ and generalized Baire space $\kappa^\kappa$ for an infinite cardinal $\kappa>\omega$ and study special subsets of these spaces. In the recent years the theory of the generalized Cantor and Baire spaces was extensively developed (see, e.g. \cite{pllmps:hdgbs}, \cite{sfthvk:gdstct}, \cite{sfgl:nii}, \cite{gl:afnrprl}, \cite{gl:gssm}, \cite{ss:pnii}, \cite{sssc:grrfic},  \cite{sf:hdsti}, \cite{sf:idst}, \cite{sf:ccu} and many others). An important part of the research in this subject is an attempt to transfer the results in set theory of the real line to $2^\kappa$ and $\kappa^\kappa$ (the list of open questions can be found in \cite{glblis:qgbs}). Despite the rapid development in this theory, the authors are not aware of any thorough research in the subject of special subsets in $2^\kappa$ or $\kappa^\kappa$. Known results are related mainly to the ideal of strongly null sets (see \cite{ah:nsgbs} and \cite{ahss:smzs}).

Throughout this paper, unless it is stated otherwise, we assume that $\kappa$ is an uncountable regular cardinal number.  

\subsubsection{Preliminaries}

We consider the space $2^\kappa$, called \df{$\kappa$-Cantor space (or the generalized Cantor space)}, endowed with so called bounded topology given by a base $\{[x]\colon x\in 2^{<\kappa}\}$, where for $x\in 2^{<\kappa}$, $[x]=\{f\in 2^\kappa\colon f\obc\dom x= x\}$.

Similarly, the space $\kappa^\kappa$ along with bounded topology given by a base $\{[x]\colon x\in \kappa^{<\kappa}\}$, where for $x\in \kappa^{<\kappa}$, $[x]=\{f\in \kappa^\kappa\colon f\obc\dom x= x\}$. is called \df{$\kappa$-Baire space (or the generalized Baire space)}. 

Throughout this paper, let $\Kont\in\{2,\kappa\}$. Therefore, $\Kont^\kappa$ denotes the generalized Cantor space or the generalized Baire space. 

If we additionally assume that $\kappa^{<\kappa}=\kappa$, then the above base has cardinality $\kappa$. This assumption proves to be very convenient when considering the generalized Cantor space and the generalized Baire space, and is assumed throughout this paper, unless stated otherwise (see e.g.  \cite{sfthvk:gdstct}).

The space $2^\kappa$ will also be treated as a~vector space over $\Z_{2}$. In particular, for $A,B\subseteq 2^{\kappa}$, let $A+B=\{t+s\colon t\in A, s\in B\}$.  Let $\bzero\in 2^\kappa$ be such that $\bzero(\alpha)=0$ for all $\alpha<\kappa$, $\bjeden\in 2^\kappa$ be such that $\bjeden(\alpha)=1$ for all $\alpha<\kappa$, and let $\QQ=\{x\in 2^\kappa\colon \exists_{\alpha<\kappa}\forall_{\alpha<\beta<\kappa} x(\beta)=0\}$. Similarly, if $x,y\in \kappa^\kappa$, then $x+y\in \kappa^\kappa$ is such that $x(\alpha)+y(\alpha)=(x+y)(\alpha)$ (in the sense of addition of ordinal numbers), for all $\alpha<\kappa$.

Notice that 
the bases defined above consist of clopen sets. Moreover, an intersection of less than $\kappa$ of basic sets is a~basic set or an empty set. Thus, an intersection of less than $\kappa$ open sets is still open. Notice also that there are $2^{\kappa}$ closed sets in those spaces.

Additionally, (under the assumption $\kappa^{<\kappa}=\kappa$) there exists a~family $\oF$ of subsets of $\kappa$ such that $|\oF|=2^\kappa$, and for all $A,B\in\oF$, $|A\cap B|<\kappa$ if $A\neq B$. Indeed, let $b\colon 2^{<\kappa}\to \kappa$ be a~bijection. Then $\oF=\left\{b\left[\{x\obc\alpha\colon\alpha<\kappa\}\right]\colon x\in 2^\kappa\right\}$
is such a family.

A~$T_1$ topological space is said to be \df{$\kappa$-additive} if for any $\alpha<\kappa$, an intersection of an $\alpha$-sequence of open subsets of this space is open. Various properties of $\kappa$-additive spaces were considered by R.~Sikorski in \cite{rs:rstshp}. The generalized Cantor  and  Baire spaces are examples of $\kappa$-additive spaces. It is also easy to see that every $\kappa$-additive topological space $X$ with clopen base of cardinality $\kappa$ is homeomorphic to a~subset of $2^\kappa$.

It is also easy to see that $A\subseteq 2^\kappa$ is closed if and only if $A=[T]$ for some tree $T\subseteq 2^{<\kappa}$. 
A similar fact is also true in the generalized Baire space. 
For a closed $A\subseteq \Kont^\kappa$, a~tree $T\subseteq \Kont^{<\kappa}$ such that $A=[T]$ is denoted by $T_A$.

The family of \df{$\kappa$-Borel sets} is the smallest family of subsets of $\Kont^\kappa$ containing all open sets and closed under complementation and under taking intersections of size $\kappa$. The family of such sets is denoted here by $\oB_\kappa$.  A function $f\colon \Kont^\kappa\to \Kont^\kappa$ is \df{$\kappa$-Borel measurable} if for any $s\in \Kont^{<\kappa}$, $f^{-1}[[s]]$ is $\kappa$-Borel.

We say that a~set is \df{$\kappa$-meagre} if it is a~union of at most $\kappa$ nowhere dense (in the bounded topology) sets. Notice that the generalization of the  Baire  category  theorem holds,  namely $2^\kappa$ is  not $\kappa$-meagre (see \cite[Theorem~xv]{rs:rstshp}), and neither is $\kappa^\kappa$. Obviously, the analogue of Kuratowski-Ulam theorem holds also for $\kappa$-meagre sets. The family of all $\kappa$-meagre sets in $2^\kappa$ or $\kappa^\kappa$ is denoted by $\MM_\kappa$ (the underlying space will be clear from the context). 

We let $\cof (\oM_\kappa)=\min\left\{|\oA|\colon \oA\subseteq \oM_\kappa\land \forall_{A\in\oM_\kappa} \exists_{B\in \oA} A\subseteq B\right\}$,
and \[\cov (\oM_\kappa)=\min\left\{|\oA|\colon \oA\subseteq \oM_\kappa\land \bigcup \oA = \Kont^\kappa\right\}.\]

Notice that if $\langle x_\alpha\rangle_{\alpha<\kappa}\in (\Kont^\kappa)^\kappa$ is a~sequence of points in $\Kont^\kappa$ such that for all $\xi<\kappa$, there exists $\delta_\xi<\kappa$ such that for all $\delta_\xi\leq\alpha,\beta<\kappa$, $x_\alpha\obc \xi=x_\beta\obc \xi$,  then there exists $x\in \Kont^\kappa$ which is a~(topological) limit of $\langle x_\alpha\rangle_{\alpha<\kappa}$ (i.e. for every open set $U$ with $x\in U$, there exists $\xi<\kappa$ such that for all $\xi<\alpha<\kappa$, $x_\alpha\in U$). Indeed, take $x=\bigcup_{\xi<\kappa}x_{\delta_\xi}\obc \xi$.

Obviously, if $C\subseteq \Kont^\kappa$ is closed, and $\langle x_\alpha\rangle_{\alpha<\kappa}\in (\Kont^\kappa)^\kappa$ is a~sequence of points of $C$ with limit $x\in \Kont^\kappa$, then $x\in C$ as well. Therefore, if $\langle C_\alpha\rangle_{\alpha\in\kappa}$ is a~sequence of non-empty closed sets such that $C_\beta\subseteq C_\alpha$, when $\alpha<\beta<\kappa$, and such that there exist increasing sequences $\langle \xi_\alpha\rangle_{\alpha\in \kappa}\in\kappa^\kappa$, $\langle s_\alpha\rangle_{\alpha\in \kappa}\in (\Kont^{<\kappa})^\kappa$, with $C_\alpha\subseteq [s_\alpha]$ and $s_\alpha\in \Kont^{\xi_\alpha}$, then there is $x\in \Kont^\kappa$ such that $\bigcap_{\alpha<\kappa} C_\alpha=\{x\}$.

Obviously, spaces $2^\kappa\times 2^\kappa$ and $2^\kappa$ are homeomorphic, and the canonical homeomorphism between them is given by the canonical well-order of $2\times\kappa$, $g\colon 2 \times\kappa\to \kappa$.

\subsubsection{Cardinal coefficients in $2^\kappa$}

The~statement $2^\kappa=\kappa^+$ is known as the \df{Continuum Hypothesis for $\kappa$} and is denoted by $CH_\kappa$.

Recall that $\diamondsuit_\kappa(E)$ for $E\subseteq \kappa$ is the following principle: there exists a~sequence $\langle S_\alpha\rangle_{\alpha\in E}$ such that $S_\alpha\subseteq \alpha$ for all $\alpha\in E$, and the set 
$\left\{\alpha\in E\colon X\cap \alpha=S_\alpha\right\}$
is a stationary subset of $\kappa$ for every $X\subseteq \kappa$ (see e.g. \cite{tj:st}[Chapter~23]). The principle $\diamondsuit_\kappa(\kappa)$ is simply denoted by $\diamondsuit_\kappa$ (and called the \df{diamond principle for $\kappa$}).

If $f,g\in \kappa^\kappa$, then we write $f\leq^\kappa g$ if there exists $\alpha<\kappa$ such that for all $\beta<\kappa$ if $\beta>\alpha$, then $f(\beta)\leq g(\beta)$. In this case we say that $f$ is \df{eventually dominated} by $g$.

Analogously as in the case of $\omega^\omega$, one can define cardinals related to the order $\leq^\kappa$.
The following two cardinals play an important role: $\bont_\kappa=\min\{|\oA|\colon \oA\subseteq \kappa^\kappa\land \lnot\exists_{f\in \kappa^\kappa}\forall_{g\in A}g\leq^\kappa f\}$,
and
$\dont_\kappa=\min\{|\oA|\colon \oA\subseteq \kappa^\kappa\land \forall_{f\in \kappa^\kappa}\exists_{g\in A}f\leq^\kappa g\}$,
which are called the \df{bounding and dominating number for $\kappa$}, respectively. Obviously, $\kappa<\bont_\kappa\leq\dont_\kappa\leq 2^\kappa$. The cardinal coefficients of an uncountable cardinal number were extensively studied in \cite{jbabsfdm:cduc}.

\subsubsection{$\kappa$-Compactness}

Not all the results of theory of the real line can be easily generalized to the case of $2^\kappa$. One of the main obstacles is the notion of compactness. We shall say that a~topological space $X$ is \df{$\kappa$-compact} (or $\kappa$-Lindel\"{o}f) if every open cover of $X$ has a~subcover of cardinality less than $\kappa$ (see \cite{dmds:asret}, \cite{hhsn:sh}). Obviously, the Cantor space $2^\omega$ is $\omega$-compact (i.e. compact in the traditional sense). But it is not always the case that $2^\kappa$ is $\kappa$-compact. Recall that a~cardinal number $\kappa$ is weakly compact if it is uncountable, and for every two-colour colouring of the set of all two-element subsets of $\kappa$, there exists a~set $H\subseteq \kappa$ of cardinality $\kappa$, which is homogeneous (every two-element subset of $H$ have the same colour in the considered colouring) (see \cite{tj:st}). Recall also that every weakly compact cardinal is strongly inaccessible. Actually, the generalized Cantor space $2^\kappa$ is $\kappa$-compact if and only if $\kappa$ is a~weakly compact cardinal (see \cite{dmds:asret}).

And there is even more to that. 
The generalized Cantor space $2^\kappa$ and the generalized Baire spaces $\kappa^\kappa$ are homeomorphic if and only if $\kappa$ is not a~weakly compact cardinal (see \cite{hhsn:sh}).

\subsubsection{Perfect sets in $\Kont^\kappa$}

A~set $P\subseteq \Kont^\kappa$ is a~\df{perfect set} if it is closed and has no isolated points. Notice that a~set $P\subseteq \Kont^\kappa$ is perfect if and only if $T_P$ is a~perfect tree.

A~perfect tree $T$ will be called \df{$\kappa$-perfect} if for every limit $\beta<\kappa$, and $t\in \Kont^\beta$ such that for all $\alpha<\beta$, $t\obc\alpha\in T$, we have $t\in T$. 

Notice that every $\kappa$-perfect tree $T\subseteq 2^{<\kappa}$ is order-isomorphic with $2^{<\kappa}$. 

A~set $P\subseteq \Kont^\kappa$ is \df{$\kappa$-perfect} if $P=[T]$ for a~$\kappa$-perfect tree $T$. Obviously, every $\kappa$-perfect set is perfect. On the other hand, the converse does not hold. 

Another major difference between $2^\kappa$ and $2^\omega$ is related to the perfect set property of an analytic set. In  $2^\omega$ every uncountable analytic set contains a~perfect set. However, the generalization of this theorem for $2^\kappa$ may not be true even for closed sets. There may even exist a perfect set which does not contain a~$\kappa$-perfect set. Recall that a~tree $T\subseteq 2^{<\kappa}$ is a~\df{$\kappa$-Kurepa} tree if:
\begin{enumerate}[(1)]
\item $|[T]|>\kappa$,
\item if $\alpha$ is uncountable, then $|T\cap 2^\alpha|\leq |\alpha|$.
\end{enumerate}
If $T$ is a~$\kappa$-Kurepa tree, then $[T]$ is an example of a~closed set of cardinality larger than $\kappa$, with no $\kappa$-perfect subsets (see e.g. \cite{glblis:qgbs, sf:hdsti}). 

Fortunately, one can see that every $\kappa$-comeagre set contains a~$\kappa$-perfect set. Indeed, if $G=\bigcup_{\alpha<\kappa} G_\alpha$ with $G_\alpha$ nowhere dense, we choose by induction $\left<t_s\right>_{s\in \Kont^{<\kappa}}$ such that $t_s \in \Kont^{<\kappa}$, and for $s,s'\in \Kont^{<\kappa}$, $s\subsetneq s'$ if and only if $t_s\subsetneq t_{s'}$. Indeed, let $t_\zbp$ be such that $[t_\zbp]\cap G_0=\zbp$. Then, given $t_s$, $s\in \Kont^\alpha$, let $t'_s\supsetneq t_s$ be such that $[t'_s]\cap G_{\alpha+1}= \zbp$. For $\xi\in\Kont$, set $t_{s^\frown \xi}=t_s^{\prime\frown}\xi$. For limit $\beta<\kappa$, and $s\in \Kont^\beta$, let $t'_s=\bigcup_{\alpha<\beta}t_{s\obc\alpha}$. Let $t_s\supsetneq t'_s$ be such that $[t_s]\cap G_\beta=\zbp$. Finally, put $T=\bigcup_{\alpha<\kappa}\{t\in \Kont^{<\kappa}\colon t\subseteq t_{s}, s\in \Kont^{\alpha}\}$. Obviously, $T$ is a~$\kappa$-perfect tree, so $P=[T]_{\kappa}$ is a~$\kappa$-perfect subset of $\Kont^\kappa\setminus G$.

In \cite{pllmps:hdgbs}, the authors consider the following properties of subsets of $\kappa^\kappa$, which generalize two well-known conditions from the classical descriptive set theory.  A subset $X\subseteq \kappa^\kappa$ satisfies \df{the Hurewicz dichotomy} if either $X$ is included in a $\kappa$-union of $\kappa$-compact sets or contains a closed (in $\kappa^\kappa$) set which is homeomorphic to $\kappa^\kappa$. We shall say that $X\subseteq \kappa^\kappa$ has \df{the perfect set property} if either $|X|\leq \kappa$ or $X$ contains a closed in $\kappa^\kappa$ set homeomorphic to $2^\kappa$.

It is proven in \cite{pllmps:hdgbs}, that for an uncountable cardinal $\kappa$, there is a partial order which forces that every $\Sigma_1^1$ set in $\kappa^\kappa$ satisfies the Hurewicz dichotomy. Moreover, the authors show  that one can find a class forcing notion (under GCH) which forces the above property for every uncountable regular $\kappa$, and such that it also forces the failure of the $\kappa$-perfect set property for closed subsets of $\kappa^\kappa$. 

On the other hand, the same paper provides a proof of the fact that after adding a Cohen subset of $\mu$ over $V$, where $\mu$ is a regular cardinal, and $\mu^{<\mu}=\mu<\kappa=\kappa^{<\kappa}$, the set $\kappa^\kappa$ in $V$ is an example of a closed set for which the Hurewicz dichotomy fails. Moreover, it is shown that in the constructible universe for every uncountable cardinal $\kappa$, there is a closed set in $\kappa^\kappa$ which does not satisfy the Hurewicz dichotomy.

\section{Special subsets of $2^\kappa$ and $\kappa^\kappa$: simple generalizations}\label{chsimple}
The aim of this section is to generalize to the case of $2^\kappa$ or $\kappa^\kappa$ certain notions of special subsets defined for $2^\omega$, and to check their properties and relations between them. Most of the results presented here have their counterparts in the standard case of $2^\omega$, and if so we give a~reference in the form ($\omega$: [n]). 

The results presented in this section consist of relatively simple generalizations of some results summarized in \cite{am:ssrl} and \cite{lb:srl} to the case of $2^\kappa$.

\subsection{Lusin sets for $\kappa$}

Recall that a~set $L\subseteq 2^{\omega}$ is a~\df{$\kappa$-Lusin set} (for $\omega<\kappa\leq 2^\omega$) if for any meagre set $X$, $|L\cap X|<\kappa$, but $|L|\geq \kappa$. An~$\aleph_1$-Lusin set is simply called a~\df{Lusin set}. This idea was introduced independently in \cite{nl:pb} and \cite{pm:tkm}. The existence of a~Lusin set is independent from ZFC. It is easy to see that under CH such a~set exists. Indeed, enumerate all closed nowhere dense sets and inductively pick a~point from a~complement of each such set distinct from all the points chosen so far. The same can be easily done if $\cov(\MM)=\cof(\MM)=\aleph_1$ (see e.g. \cite{lb:srl}).

Let $\kappa<\lambda\leq 2^\kappa$. A~set $L\subseteq \Kont^\kappa$ such that $|L|\geq \lambda$, and if $X\subseteq \Kont^\kappa$ is any $\kappa$-meager set, then $|X\cap L|< \lambda$ will be called a~\df{$\lambda$-$\kappa$-Lusin set}. A~$\kappa^+$-$\kappa$-Lusin set is simply called a~\df{Lusin set for $\kappa$}.

\begin{theo}[$\omega$: \cite{lb:srl}, Theorem 8.26]\label{p-lusin}
If $\lambda=\cov(\MM_\kappa)=\cof(\MM_\kappa)$, then there exists a~$\lambda$-$\kappa$-Lusin set $L\subseteq 2^\kappa$.
\end{theo}

Proof: 
The proof is a straightforward generalization of the proof in the case $\kappa=\omega$, which can be found in \cite[Theorem 8.26]{lb:srl}.\ \hfill $\square$

Obviously, since $2^\kappa\subseteq \kappa^\kappa$, we get that under the above conditions there exists a exists a~$\lambda$-$\kappa$-Lusin set $L\subseteq \kappa^\kappa$. Also, immediately we get the following corollary.

\begin{cor}
Assume $CH_\kappa$. Then there exists a~Lusin set for $\kappa$ in $\Kont^\kappa$.
\end{cor}

Proof: Clear.\ \hfill $\square$

On the other hand, the existence of a~$\lambda$-$\kappa$-Lusin set constrains the value of $\cov(\MM_\kappa)$.

\begin{prop}[$\omega$: \cite{lb:srl}, Theorem 8.35]
Assume that $\lambda$ is a~regular cardinal and $\kappa<\lambda\leq 2^\kappa$. If $L$ is a~$\lambda$-$\kappa$-Lusin set in $\Kont^\kappa$, then $|L|\leq \cov(\MM_\kappa)$.
\end{prop}

Proof: 
The proof is a straightforward generalization of the proof in the case $\kappa=\omega$, which can be found in \cite[Theorem 8.35]{lb:srl}.\ \hfill $\square$

\begin{cor}
Assume that $\lambda$ is a~regular cardinal, $\kappa<\lambda\leq 2^\kappa$, and that there exists a~$\lambda$-$\kappa$-Lusin set $L$. Then $\non(\MM_\kappa)\leq \lambda\leq \cov(\MM_\kappa)$.
\end{cor}
Proof: Clear.\ \hfill $\square$

\begin{cor}
Assume that $\lambda$ is a~regular cardinal, $\kappa<\lambda\leq 2^\kappa$. There exists a model of ZFC such that there are no~$\lambda$-$\kappa$-Lusin sets.
\end{cor}

Proof: \cite[Theorem~49]{jbabsfdm:cduc} states that there is a model of ZFC in which $\cov(\MM_\kappa)=\kappa^+$ and $\non(\MM_\kappa)=\kappa^{+++}$.\ \hfill $\square$

\subsection{$\kappa$-strongly measure zero in $\Kont^\kappa$}\label{s-sn}

In the classical set theory of real line a~set $A$ is called \df{strongly null} (strongly of measure zero) if for any sequence of positive $\varepsilon_{n}>0$, there exists a~sequence of open sets $\left<A_{n}\right>_{n\in\omega}$, with $\diam A_{n}<\vare_{n}$ for $n\in\omega$, and such that $A\subseteq\bigcup_{n\in\omega}A_{n}$. We denote the class of such sets by $\SN$. The idea of strongly null was introduced for the first time in \cite{eb:cemn}, and \df{Borel conjectured} that all $\SN$ sets are countable. This hypothesis turned out to be independent from ZFC (see \cite{rl:cbc}). It is easy to see that a~set $A$ is strongly null if and only if for any sequence of positive $\varepsilon_{n}>0$, there exists a~sequence of open sets $\left<A_{n}\right>_{n\in\omega}$, with $\diam A_{n}<\vare_{n}$ for $n\in\omega$, and such that 
$A\subseteq\bigcap_{m\in\omega}\bigcup_{n>m}A_{n}$.

Galvin, Mycielski and Solovay (in \cite{fgjmrs:smzs}) proved that a~set $A\in\SN$ (in $2^{\omega}$) if and only if for any meagre set $B$, there exists $t\in 2^{\omega}$ such that $A\cap(B+t)=\zbp$. 

A~set $A\subseteq \Kont^\kappa$ will be called \df{$\kappa$-strongly measure zero} ($\kappa$-strongly null, $\SN_\kappa$) if for every $\langle \xi_\alpha \rangle_{\alpha<\kappa}\in \kappa^\kappa$, there exists $\left<x_\alpha\right>_{\alpha< \kappa}$ such that $x_\alpha\in \Kont^{\xi_\alpha}$, $\alpha<\kappa$ and $A\subseteq \bigcup_{\alpha< \kappa} [x_\alpha]$ (see also \cite{ah:nsgbs} and \cite{ahss:smzs}). Obviously if $A\in [\Kont^\kappa]^{\leq\kappa}$, then $A\in\SN_\kappa$. 

The following well-known characterization of strongly null sets can be generalized to $\Kont^\kappa$.

\begin{prop}[$\omega$: \cite{lb:srl}]
If $A\in\SN_\kappa$, and $\langle \xi_\alpha \rangle_{\alpha<\kappa}\in \kappa^\kappa$, there exists $\left<x_\alpha\right>_{\alpha< \kappa}\in(\Kont^\kappa)^\kappa$ such that $x_\alpha\in \Kont^{\xi_\alpha}$ for all $\alpha<\kappa$, and
$A\subseteq \bigcap_{\alpha<\kappa}\bigcup_{\alpha<\beta< \kappa} [x_\beta]$.
\end{prop}

Proof: The proof is a straightforward generalization of the proof in the case $\kappa=\omega$, which can be found in \cite{lb:srl}.\ \hfill $\square$

In particular, the family of $\SN_\kappa$ sets forms a~$\kappa^+$-complete ideal.

Notice also that, obviously, $\Kont^\kappa\notin \SN_\kappa$. 

The \df{Generalized Borel Conjecture for $\kappa$ (GBC($\kappa$))} states that $\SN_\kappa=[\Kont^\kappa]^{\leq\kappa}$. 

Some properties of $\SN_\kappa$ of sets were considered in \cite{ahss:smzs}. In particular, it is proven there that if $\kappa$ is a~successor cardinal, then $\SN_\kappa$ is a~$\bont_\kappa$-additive ideal. Under Generalized Martin Axiom for $\kappa$ (GMA($\kappa$), see \cite{ss:wgmahc}), $\bont_\kappa=2^\kappa$, so then $\SN_\kappa$ is $2^\kappa$-additive. Finally in \cite{ss:wgmahc}, it is proven that GBC($\kappa$) fails for a successor $\kappa$.

We study some other properties of $\kappa$-strongly measure zero sets.

\begin{prop}[$\omega_1$: \cite{ah:nsgbs}, Proposition 9.7]
The family of all closed subsets of $\Kont^\kappa$ which are not $\SN_\kappa$ does not satisfy $2^\kappa$-chain condition.
\end{prop}

Proof: Since $\kappa^{<\kappa}=\kappa$, we can consider a~family $\oF$ of subsets of $\kappa$ such that $|\oF|=2^\kappa$, and for all $X,Y\in\oF$, $|X\cap Y|<\kappa$ if $X\neq Y$. The rest of the proof is a straightforward generalization of the proof of \cite[Proposition 9.7]{ah:nsgbs}. 
\hfill $\square$

\begin{prop}[$\omega$: \cite{am:ssrl}, Theorem~8.5]
Assume $CH_\kappa$. Then there exists a~Lusin set for $\kappa$ $L\subseteq 2^\kappa$ such that $L\times L\notin \SN_\kappa$.
\end{prop}

Proof: Let $\left<X_\alpha\colon \alpha<\kappa^+\right>$ be an enumeration of all closed nowhere dense sets, and let $\{y_{\alpha}\colon \alpha<\kappa^+\}=2^{\kappa}$. Inductively, for $\alpha<\kappa^+$, choose \[x_\alpha, x'_{\alpha}\in 2^\kappa\setminus \left(\{x_\beta\colon\beta<\alpha\}\cup\{x'_\beta\colon\beta<\alpha\}\cup \bigcup_{\beta<\alpha}X_\beta\right)\]
such that $x_\alpha+x'_\alpha=y_\alpha$. This is possible, since $F_{\alpha}=\{x_\beta\colon\beta<\alpha\}\cup\{x'_\beta\colon\beta<\alpha\}\cup \bigcup_{\beta<\alpha}X_\beta$ is $\kappa$-meagre, so $(y_{\alpha}+ F_{\alpha})\cup F_{\alpha}$ is also $\kappa$-meagre. Thus, there exists $x_{\alpha}\notin (y_{\alpha}+ F_{\alpha})\cup F_{\alpha}$. Let $x'_{\alpha}=x_{\alpha}+y_{\alpha}$. Then also $x'_{\alpha}\notin F_{\alpha}$.

Obviously, $L=\{x_\alpha\colon\alpha<\kappa\}\cup\{x'_\alpha\colon\alpha<\kappa\}$ is a~Lusin set for $\kappa$. Nevertheless, $L\times L$ is not a~$\SN_\kappa$ set. Indeed, let $f\colon 2^{\kappa}\times 2^{\kappa}\to 2^\kappa$ be given by $f(x,x')=x+x'$. Notice that if $\alpha<\kappa$ is a~limit ordinal, then $g(0,\alpha)=\alpha$, where $g$ is the canonical well-ordering of $2\times \kappa$. Therefore, if $x\in 2^\beta$, for $\omega\leq \beta<\kappa$, then $[x]$ when considered as a~subset of $2^\kappa\times 2^\kappa$ (we identify $2^\kappa\times 2^\kappa$ and $2^\kappa$ by applying $g$) is contained in $[y]\times [z]$, where $y,z\in 2^\alpha$ with $\alpha$ a~limit ordinal such that $\alpha\leq \beta <\alpha+\omega$. This implies that $f[[x]]\subseteq [y+z]$, and  thus if $X\subseteq 2^\kappa\times 2^\kappa$ is $\kappa$-strongly null, then $f[X]$ is $\SN_\kappa$ as well. But $f[L\times L]=2^\kappa$, so $L\times L\notin \SN_\kappa$. \ \hfill $\square$ 

Next, we study the possibility of a generalization of Galvin, Mycielski and Solovay \cite{fgjmrs:smzs} characterization of strongly null sets. One of the implications can be proved under no additional assumptions. Before finalization of this paper, the authors became aware that results of Proposition~\ref{gms1}, Lemma~\ref{lem-gms}, Theorem~\ref{gms} were independently proved by W.~Wohofsky, and can be found also in~\cite{ww:ssrnvbc}. Nevertheless, we present here this result with a proof for the sake of completeness.

\begin{prop}[$\omega$: \cite{am:ssrl}, Theorem 3.5]\label{gms1}
Let $A\subseteq 2^\kappa$ be such that for any nowhere dense set $F$, there exists $x\in 2^\kappa$ such that $(x+A)\cap F=\zbp$. Then, $A$ is $\SN_\kappa$.
\end{prop}

Proof: The proof is a straightforward generalization of the proof in the case $\kappa=\omega$, which can be found as a part of the proof of \cite[Theorem 3.5]{am:ssrl}.\ \hfill $\square$

The reversed implication can be proven under the assumption that $\kappa$ is a~weakly compact cardinal.

\begin{lemma}[$\omega$: \cite{am:ssrl}, Lemma 3.5.1]\label{lem-gms}
Assume that $\kappa$ is weakly compact. For any closed nowhere dense set $C\subseteq 2^\kappa$ and $s\in 2^{<\kappa}$, there exists $\xi<\kappa$ and $F\subseteq \{s'\in 2^{<\kappa}\colon s\subsetneq s'\}$ with $|F|<\kappa$ such that for any $t\in 2^\xi$, there exists $s'\in F$ such that $([s']+[t])\cap C=\zbp$.
\end{lemma}

Proof: The proof is a relatively easy generalization of the proof of \cite[Lemma~3.5.1]{am:ssrl}.  
\ \hfill $\square$

\begin{theo}[$\omega$: \cite{am:ssrl}, Theorem 3.5]\label{gms}
Assume that $\kappa$ is a~weakly compact cardinal, and $A\subseteq 2^\kappa$ is $\SN_\kappa$. Then for any $\kappa$-meagre set $F$, there exists $x\in 2^\kappa$ such that $(x+A)\cap F=\zbp$.
\end{theo}

Proof: Let $F=\bigcup_{\alpha<\kappa} C_\alpha$ with $C_\alpha$ closed nowhere dense sets. We can assume that $C_{\alpha}\subseteq C_\beta$ if $\alpha<\beta$. 

We construct inductively a~tree $T\subseteq \kappa^{<\kappa}$, along with sequences $\langle \delta_u\rangle_{u\in T},\langle \xi_u\rangle_{u\in T}\in \kappa^T$, and $\left<s_u\right>_{u\in T}\in \left(2^{<\kappa}\right)^T$  such that:
\begin{enumerate}[(a)]
\item if $u\in T\cap\kappa^\beta$, $\beta<\kappa$, then $\{u'\in T\cap\kappa^{\beta+1}\colon u\subseteq u'\}=\{u^\frown\alpha\colon \alpha<\delta_u\}$,
\item for any $u,u'\in T$ if $u\subsetneq u'$, then $s_u\subsetneq s_{u'}$,
\item for any $u\in T\cap \kappa^\beta$, $\beta<\kappa$, and $t\in 2^{\xi_u}$, there exists $\alpha<\delta_u$ such that $\left([s_{u^\frown \alpha}]+[t]\right)\cap C_{\beta}=\zbp$. 
\end{enumerate}
Precisely, let $s_\zbp=\zbp$. If $u\in T\cap \kappa^\beta$, $\beta<\kappa$, apply Lemma~\ref{lem-gms} to $C_\beta$ and $s_u$ to get $\xi_u<\kappa$, $F_u\subseteq \{s'\in 2^{<\kappa}\colon s\subseteq s'\}$ with $|F_u|=\delta_{u}<\kappa$ such that for any $t\in 2^{\xi_u}$, there exists $s'\in F_u$, so that $([s']+[t])\cap C_\beta=\zbp$. Fix an enumeration $F_u=\{s'_{u,\alpha}\colon \alpha<\delta_u\}$, and put $\{u'\in T\cap\kappa^{\beta+1}\colon u\subseteq u'\}=\{u^\frown\alpha\colon \alpha<\delta_u\}$. For all $\alpha<\delta_u$, let $s_{u^\frown\alpha}=s'_{u,\alpha}$. If $\beta<\kappa$ is a~limit ordinal, let 
$T\cap \kappa^\beta=\{u\in \kappa^\beta\colon \forall_{\alpha<\beta} u\obc\alpha \in T\}$.
Also, for $u\in T\cap \kappa^\beta$, let $s_u=\bigcup_{\alpha<\beta} s_{u\obc \alpha}$. 

Next, define $\langle \delta_\alpha\rangle_{\alpha<\kappa},\langle \xi_\alpha\rangle_{\alpha<\kappa}$ in the following way. For $\alpha<\kappa$, let $\delta_\alpha=\bigcup_{u\in T\cap \kappa^\alpha}\delta_u,$ and $\xi_\alpha=\bigcup_{u\in T\cap \kappa^\alpha}\xi_u$. Notice that for all $\alpha<\kappa$, $\delta_\alpha, \xi_\alpha<\kappa$. Indeed, if it holds for $\alpha<\kappa$, then $|T\cap\kappa^{\alpha+1}|=\delta_\alpha<\kappa$, so $\delta_{\alpha+1}, \xi_{\alpha+1}<\kappa$, since $\kappa$ is regular. If $\alpha$ is a~limit ordinal, then $T\cap \kappa^\alpha\subseteq \delta^\alpha$ with $\delta=\bigcup_{\beta<\alpha}\delta_{\beta}<\kappa$. And $\delta^\alpha<\kappa$, as $\kappa$ is strongly inaccessible (every weakly compact cardinal is strongly inaccessible).

Assuming that $A$ is $\SN_\kappa$, there exists $\left<x_\alpha\right>_{\alpha\in \kappa}$ such that $x_\alpha\in 2^{\xi_\alpha}$, $\alpha\in \kappa$ and 
$A\subseteq \bigcap_{\beta<\kappa}\bigcup_{\beta<\alpha<\kappa} [x_\alpha]$.
By induction, construct $y\in\kappa^\kappa$ such that: 
\begin{enumerate}[(a)]
\item for all $\alpha<\kappa$, $y\obc \alpha \in T$,
\item for all $\alpha<\kappa$, $\left([s_{y\obc(\alpha+1)}]+[x_\alpha]\right)\cap C_\alpha=\zbp$.
\end{enumerate}
Precisely, let $y(\alpha)<\delta_{y\obc\alpha}$ be such that 
$\left([s_{y\obc\alpha^\frown y(\alpha)}]+[x_\alpha]\right)\cap C_\alpha=\zbp$.
Notice that if $\alpha$ is a~limit ordinal, then $y\obc\alpha=\bigcup_{\beta<\alpha} y\obc\beta\in T$.

Finally, let 
$x=\bigcup_{\alpha<\kappa} s_{y\obc\alpha}\in 2^\kappa$.
Notice that for all $\beta\leq\alpha<\kappa$, we have $(x+[x_\alpha])\cap C_\beta=\zbp$. Therefore, 
$(x+A)\cap F=\zbp$.
\ \hfill $\square$

We do not know whether for other cardinals $\kappa$ the statement of this theorem can be generalized or refuted.

\begin{prob}\label{probsn}
Does the above theorem hold for a cardinal $\kappa$ which is not weakly compact? Which assumptions on $\kappa$ refute the hypothesis of this theorem?
\end{prob}

A related open question, in the case of $\kappa=\omega_1$ was  stated in \cite{ah:nsgbs}[Conjecture~7.9]. In \cite{ah:nsgbs}, the author conjectures that the class of $\kappa$-strongly measure zero sets coincides with the family of stationary strongly measure zero sets, which, as he proves in \cite{ah:nsgbs}[Theorem~7.8], under some additional assumptions, satisfy the hypothesis of Galvin-Mycielski-Soloway Theorem. A set $A$ is defined to be stationary strongly measure zero, if for every $\langle \xi_\alpha \rangle_{\alpha<\kappa}\in \kappa^\kappa$, there exists $\left<x_\alpha\right>_{\alpha< \kappa}$ such that $x_\alpha\in \Kont^{\xi_\alpha}$, $\alpha<\kappa$ and for every club set $C\subseteq \kappa$,  $A\subseteq \bigcup_{\alpha\in C} [x_\alpha]$.

On the other hand, in \cite{ww:ssrnvbc}[Question~5.20], the same open problem, as in Question~\ref{probsn} is formulated. In \cite{ww:ssrnvbc}, the author conjectures, with regard to the importance of the assumption of compactness in the proofs of analogues of Galvin-Mycielski-Soloway Theorem, that for regular uncountable cardinals which are not weakly compact,  Theorem~\ref{gms} does not hold. 

Obviously, an affirmative answer on Problem~\ref{probsn} for $\kappa$ would follow from $GBC(\kappa)$. But, consistency of $GBC(\kappa)$ is also an open problem (see \cite{ahss:smzs}).

The above theorem and propositions imply the following corollaries (see \cite[Corollary 8.14]{lb:srl}).

\begin{cor}[$\omega$: \cite{lb:srl}]
Assume that $\kappa$ is weakly compact, and $A,B\subseteq 2^\kappa$ are such that $|A|<\add(\MM_\kappa)$ and $B\in \SN_\kappa$. Then $A\cup B\in \SN_\kappa$.
\end{cor}

Proof: The proof is similar to the proof of the \cite[Corollary 8.14 a)]{lb:srl}. Since we have to use Theorem~\ref{gms}, we assume that $\kappa$ is weakly compact.\hfill $\square$

\begin{cor}[$\omega$: \cite{lb:srl}]\label{covsn}
If $A\subseteq 2^\kappa$, and $|A|<\cov(\MM_\kappa)$, then $A\in\SN_\kappa$.
\end{cor} 

Proof: The proof is similar to the proof of the \cite[Corollary 8.14 b)]{lb:srl}. Since we only need Proposition~\ref{gms1}, no additional assumptions are needed.  \ \hfill $\square$

\subsection{$\kappa^+$-Concentrated sets}

Recall that a set $A$ is said to be \df{concentrated on a~set $C$} if $A\setminus U$ is countable for every open $U$ with $C\subseteq U$ (\cite{fr:epcdpl}). Notice that every set concentrated on a~countable set is $\SN$.

Thus, we shall say that a~set $A\subseteq \Kont^\kappa$ is \df{$\lambda$-concentrated on a~set $B\subseteq \Kont^\kappa$} (for $\kappa<\lambda\leq 2^\kappa$) if for any open set $G$ such that $B\subseteq G$, we have $|A\setminus G|<\lambda$.

The classical relation between concentrated sets, Lusin sets and strongly null sets can be easily generalized to $\kappa$.

\begin{prop}[$\omega$: \cite{am:ssrl}, Theorem 3.1(i)]\label{lusin-con}
A~set $A\subseteq \Kont^\kappa$ is a~Lusin set for $\kappa$ if and only if $|A|>\kappa$ and is $\kappa^+$-concentrated on every dense set $D\subseteq \Kont^\kappa$ with $|D|=\kappa$.
\end{prop}

Proof: The proof is similar to the proof in the case $\kappa=\omega$, which can be found in \linebreak \cite[Theorem 3.1(i)]{am:ssrl}.  \hfill $\square$

\begin{prop}[$\omega$: \cite{am:ssrl}, Theorem 3.1(ii)]\label{consn}
If a~set $A\subseteq \Kont^\kappa$ is $\kappa^+$-concentrated on a~set $B$ such that $|B|\leq \kappa$, then $A\in\SN_\kappa$.
\end{prop} 
 Proof: The proof is a straightforward generalization of the proof in the case $\kappa=\omega$ (see \cite[Theorem 3.1(ii)]{am:ssrl}).\ \hfill $\square$

\begin{cor}[$\omega$: \cite{am:ssrl, lb:srl}]
Every Lusin set for $\kappa$ in $\Kont^\kappa$ is $\SN_\kappa$
\end{cor}
Proof: Clear.\ \hfill $\square$

On the other hand, we get the following.

\begin{prop}[$\omega$: \cite{am:ssrl}, Theorem 3.4]
Assume $CH_\kappa$. Then there exists a~set $A\subseteq 2^\kappa$ such that $A\in\SN_\kappa$, but $A$ is not $\kappa^+$-concentrated on any $B\subseteq 2^\kappa$ with $|B|\leq \kappa$.
\end{prop}

Proof: The proof can be obtained by generalizing the proof in the case $\kappa=\omega$ (see \cite[Theorem 3.4]{am:ssrl}) and by noticing that every $\kappa$-comeagre set contains a $\kappa$-perfect set.  \hfill $\square$

\begin{prop}[$\omega$: \cite{lb:srl}, Theorem 8.18]
If $A\subseteq 2^\kappa$ is $\cov(\MM_\kappa)$-concentrated on an~$\SN_\kappa$ set, then $A$ is also $\SN_\kappa$.
\end{prop}

Proof: The proof is an easy generalization of the proof in the case $\kappa=\omega$, which can be found in \cite[Theorem 8.18]{lb:srl}.\ \hfill $\square$

\subsection{Perfectly $\kappa$-meagre sets and $\kappa$-$\lambda$-sets}\label{s-pn}
Among classes of special subsets of the real line, the class of perfectly meager sets plays an important role. A~set is \df{perfectly meager} if it is meager relative to any perfect set, and we denote it by $\PM$ (this concept first appeared in \cite{nl:pb}). Recall also that a~set $A$ such that every countable $B\subseteq A$ is a~relative $G_\delta$-set (\cite{kk:fes}) is called \df{$\lambda$-set}. Every $\lambda$-set is perfectly meagre. A set $A$ is \df{$\lambda'$-set}, if for every countable $B$, $A\cup B$ is a~$\lambda$-set. Obviously, every $\lambda'$-set is a~$\lambda$-set.

A~set $A\subseteq \Kont^\kappa$ is a~\df{$\kappa$-$\lambda$-set} if for any $B\subseteq A$ with $|B|\leq \kappa$ there exists a sequence $\left<B_\alpha\right>_{\alpha<\kappa}$, where $B_\alpha\subseteq \Kont^{\kappa}$ are open, and $\bigcap_{\alpha<\kappa} B_\alpha \cap A=B$. 

Furthermore, a~set $A\subseteq \Kont^\kappa$ will be called \df{perfectly $\kappa$-meagre} ($\PM_\kappa$) if for every perfect $P\subseteq \Kont^\kappa$, $A\cap P$ is $\kappa$-meagre relatively to $P$. Additionally, a~set $A\subseteq \Kont^\kappa$ will be called \df{$\kappa$-perfectly $\kappa$-meagre} ($\PkM_\kappa$) if for every $\kappa$-perfect $P\subseteq \Kont^\kappa$, $A\cap P$ is $\kappa$-meagre relatively to $P$. Obviously, if $A\in \PM_\kappa$, then $A\in \PkM_\kappa$.

\begin{prop}[$\omega$: \cite{am:ssrl}, Theorem~5.2]\label{lambdapm}
Every $\kappa$-$\lambda$-set
$A\subseteq \Kont^\kappa$ is perfectly $\kappa$-meagre.
\end{prop}

Proof:   The proof is a generalization of the proof in the case $\kappa=\omega$ (see \cite[Theorem 5.2]{am:ssrl}).\ \hfill $\square$

\begin{lemma}[$\omega$: \cite{lb:srl}, Theorem 8.10]\label{fpkmk}
Let $A\subseteq 2^\kappa$ be such that for every $\kappa$-Borel measurable function $f\colon 2^\kappa\to 2^\kappa$, $f^{-1}[A]$ is $\kappa$-meagre. Then $A$ is $\kappa$-perfectly $\kappa$-meagre. 
\end{lemma}

Proof: The proof is a generalization of the proof in the case $\kappa=\omega$ (see \cite[Theorem 8.10]{lb:srl}). Notice that we consider only $\kappa$-perfect sets.\ \hfill $\square$ 
 
The above lemma provides the proof of an analogue of Grzegorek Theorem for $\PkM_\kappa$ sets.

\begin{theo}[$\omega$: \cite{lb:srl}, Theorem~8.9]
There exists a $\PkM_\kappa$ set of cardinality $\non(\MM_\kappa)$.
\end{theo}

Proof: The proof is a relatively easy generalization of the proof, which can be found in \cite[Theorem~8.9]{lb:srl}.
\hfill $\square$

However, the case of $\PM_\kappa$ sets is less clear. In the classical case if $A=\bigcup_{\xi<\omega_1} A_\xi$  is a co-analytic but not analytic set with $A_\xi$ Borel for every $\xi$, a set $E$ such that $|E\cap A_\xi|=1$ for all $\xi$ is an example of a perfectly meagre set of cardinality $\aleph_1$. Indeed, if $P$ is a perfect set, then $A\cap P$ has the Baire property, which allows to prove that $E$ is meagre in $P$ (see \cite[Theorem~8.3]{lb:srl}). In our case, since there is always a $\kappa$-analytic subset of $2^\kappa$ without $\kappa$-Baire property (see e.g. \cite{sf:idst}, \cite{ahss:smzs}), it is not clear whether there always exists a $\PM_\kappa$ set of cardinality greater then $\kappa$.

\begin{prob}\label{pnkprob}
Is there a set $A\subseteq \Kont^\kappa$ such that $|A|=\kappa^+$ and $A\in\PM_\kappa$ in every model of ZFC?
\end{prob}

Another natural question is the following.
\begin{prob}
Does there exist a set $A\in\PkM_\kappa\setminus \PM_\kappa$?
\end{prob}

In particular, if $T$ is a $\kappa$-Kurepa tree such that $|[T]|<\non(\MM_\kappa)$, then $[T]$ is an example of such a set.

A~set $A$ will be called a~\df{$\kappa$-$\lambda'$-set} if for any $F$ such that $|F|\leq \kappa$, $A\cup F$ is a~$\kappa$-$\lambda$-set.

\begin{prop}[$\omega$: \cite{am:ssrl}, Theorem 7.2]
Every union of $\kappa$ many $\kappa$-$\lambda'$-sets is a~$\kappa$-$\lambda'$-set.
\end{prop}

Proof: 
The proof is similar to the proof in the case $\kappa=\omega$ (see \cite[Theorem 7.2]{am:ssrl}).\ \hfill $\square$ 

\begin{prop}[$\omega$: \cite{am:ssrl}, Theorem 8.2]
If $X,Y\subseteq 2^\kappa$ are $\kappa$-$\lambda$ sets, then $X\times Y$ is also a~$\kappa$-$\lambda$ set.
\end{prop}

Proof: 
The proof is similar to the proof in the case $\kappa=\omega$ (see \cite[Theorem 8.2]{am:ssrl}).\ \hfill $\square$

The above proposition can be proven analogously for $\kappa$-$\lambda'$ sets.

Recall that a~set $A\subseteq 2^\omega$ is a \df{$s_0$-set} if for any perfect set $P$ there exists a~perfect set $Q\subseteq P$ with $A\cap Q=\zbp$ (\cite{em:cfscce}). Similarly, a~set $A\subseteq \Kont^\kappa$ is a~\df{$\kappa$-$s_0$-$\kappa$-set} if for any $\kappa$-perfect set $P\subseteq \Kont^\kappa$, there exists a~$\kappa$-perfect set $Q\subseteq P$ such that $Q\cap A=\zbp$.

\begin{prop}[$\omega$: \cite{am:ssrl}, Theorem~5.9]
Every $\PkM_\kappa$ subset of $2^\kappa$ is a~$\kappa$-$s_0$-$\kappa$-set.
\end{prop}

Proof: The proof is a generalization of the proof in the case $\kappa=\omega$ (see \cite[Theorem 5.9]{am:ssrl}). \ \hfill $\square$

Similar proposition can be proven for $\PM_\kappa$ sets, but using a more careful reasoning.

\begin{prop}[$\omega$: \cite{am:ssrl}, Theorem~5.9]
Every $\PM_\kappa$ subset of $2^\kappa$ is an $s_0$-set.
\end{prop}

Proof:  If $G=\bigcup_{\alpha<\kappa} G_\alpha\subseteq P$ with $G_\alpha$ nowhere dense in $P$, where $P$ is perfect, we construct by induction a~partial function $F\colon 2^{<\kappa}\to T_P$ such that for $s,s'\in \dom F$, $s\subsetneq s'$ if and only if $F(s)\subsetneq F(s')$. Indeed, let $F(\zbp)$ be such that $[F(\zbp)]\cap G_0=\zbp$. Then, given $F(s)$, $s\in 2^\alpha\cap \dom F$, let $t_s\supsetneq F(s)$ be such that $[t_s]\cap G_{\alpha+1}= \zbp$ and $t_s\in \Spl(T)$. Set $F(s^\frown 0)=t_s\,^{\frown}0$ and $F(s^\frown 1)=t_s\,^{\frown}1$. For limit $\beta<\kappa$, and $s\in 2^\beta$ such that $s\obc\alpha\in \dom F$ for all $\alpha<\beta$, let $t_s=\bigcup_{\alpha<\beta} F(s\obc\alpha)$. If $t_s\in T_P$, then let $F(s)\supsetneq t_s$ be such that $F(s)\cap G_\beta=\zbp$. Otherwise, $s\notin \dom F$. Notice that since $G_\alpha$ is nowhere dense for all $\alpha<\kappa$, for any $s\in 2^{<\beta}\cap \dom F$ there exists $s'\in 2^\beta\cap \dom F$ such that $s\subseteq s'$. Finally, let $T_Q=\{t\in 2^{<\kappa}\colon t\subseteq F(s), s\in \dom F\}$. Obviously, $T_Q\subseteq T_P$ is a~perfect tree, so $Q=[T_Q]$ is a~perfect subset of $P\setminus G$. \ \hfill $\square$

Notice that a~set having only $\kappa$-meagre homeomorphic images may not be perfectly $\kappa$-meagre.

\begin{prop}[$\omega$: \cite{am:ssrl}, Theorem~9.2]
There exists a~set $A\subseteq 2^\kappa$ which is not $\PkM_\kappa$, but its every homeomorphic image is $\kappa$-meagre.
\end{prop}

Proof:  The proof is analogous to the proof for $\kappa=\omega$, which can be found in \cite[Theorem~9.2]{am:ssrl}.
 \ \hfill $\square$

On the other hand, for $\kappa$-$\lambda$-sets we get the following.
 
\begin{prop}[$\omega$: \cite{am:ssrl}, Lemma 9.3.1 b)]
Let $A,B\subseteq \Kont^\kappa$, and assume that $f\colon A\to B$ is a~one-to-one continuous map. If $B$ is a~$\kappa$-$\lambda$ set, then $A$ is also a $\kappa$-$\lambda$-set.
\end{prop}

Proof: The proof is a straightforward generalization of the proof for $\kappa=\omega$ (see \cite[Lemma 9.3.1 b)]{am:ssrl}).

A~similar statement can be proven for $\kappa$-$\lambda'$-sets.

\begin{prop}[$\omega$: \cite{am:ssrl}, Lemma 9.3.1 c)]
Let $X,Y\subseteq \Kont^\kappa$, and assume that $f\colon X\to Y$ is a~continuous map. Let $A\subseteq X$ and $B\subseteq Y$ be such that $B$ is a~$\kappa$-$\lambda'$-set, and $f\obc A$ is one-to-one onto $B$. Then $A$ is also a~$\kappa$-$\lambda'$ set.
\end{prop}

Proof: The proof is a generalization of the proof for $\kappa=\omega$, which can be found in \cite[Lemma 9.3.1 c)]{am:ssrl}.
\ \hfill $\square$

\subsection{$\kappa$-$\sigma$-Sets}

A set $A$ is called a \df{$\sigma$-set}, if for every $F_\sigma$ set $F$, there exists a $G_\delta$ set $G$ such that $A\cap F=A\cap G$.

Thus, a~set $A\subseteq \Kont^\kappa$ will be called \df{$\kappa$-$\sigma$-set} if for any sequence of closed sets $\langle F_\alpha\rangle_{\alpha< \kappa}$, there exists a~sequence of open sets $\langle G_\alpha\rangle_{\alpha< \kappa}$ such that
$A\cap \bigcup_{\alpha<\kappa} F_{\alpha}=A\cap \bigcap_{\alpha<\kappa} G_\alpha$.

\begin{prop}[$\omega$: \cite{lb:srl}, Theorem 8.43]\label{sigmapm}
Every $\kappa$-$\sigma$-set is $\PM_\kappa$.
\end{prop}

Proof: The proof is a straightforward generalization of the proof in the case $\kappa=\omega$, which can be found in \cite[Theorem 8.43]{lb:srl}.\ \hfill $\square$

\subsection{Cover selection principles in $\Kont^\kappa$}

A cover $\oU$ of a~topological space $A$ is \df{proper} if $A\notin \oU$. From now on we assume that all considered covers are proper.

A cover $\oU$ of a~set $A$ such that for any $C\in [A]^{<\omega}$ there exists $U\in\oU$ such that $C\subseteq U$, is called an~\df{$\omega$-cover}, and we call it a~\df{$\gamma$-cover} if 
$A\subseteq \bigcup_{n\in\omega}\bigcap_{m\geq n} U_m$.
 
The family of all $\omega$-covers (respectively, $\gamma$-covers) of $A$ is denoted by $\Omega(A)$ (respectively, $\Gamma(A)$). The family of all open covers of $A$ is denoted by $\oO(A)$. The underlying set is often omitted in this notation if it is clear from the context. 

If $\oA$ and $\oB$ are families of covers of a~topological space $X$, then $X$ satisfies the \df{$S_1(\oA,\oB)$ principle} if for every sequence $\langle \oU_n\rangle_{n\in\omega}\in\oA^\omega$, there exists $\oU=\{ U_n\colon n\in\omega\}$ with $U_n\in\oU_n,$ for all $n\in\omega$ such that $\oU\in\oB$. $X$ satisfies the \df{$U_{<\omega}(\oA,\oB)$ principle} if for every sequence $\langle \oU_n\rangle_{n\in\omega}\in\oA^\omega$ such that for every $n\in\omega$ and each $\oW\subseteq \oU_n$ finite, $\oW$ is not a~cover, there exists $\langle \oW_n\rangle_{n\in\omega}$ such that $\oW_n\in [\oU_n]^{<\omega}$, and $\{\bigcup \oW_n\colon n\in\omega\}\in\oB$. The covering principles were first systematically studied in \cite{ms:cocrt}.

A~set $A$ is a $\gamma$-set if for every open $\omega$-cover $\oU$, there exists $\oV\subseteq\oU$ which is a~$\gamma$-cover (\cite{jgzn:spc}), 
It can be proven that a~set $X$ is a~$\gamma$-set if and only if $X$ satisfies $S_1(\Omega,\Gamma)$ \cite{jgzn:spc}.

A~set $X$ is said to have the \df{Menger property} \cite{km:up} if it satisfies $U_{<\omega}(\oO,\oO)$. It has the \df{Hurewicz property} \cite{wh:fsf} if it satisfies $U_{<\omega}(\oO,\Gamma)$. Finally, it has the \df{Rothberger property} \cite{fr:vec} if it satisfies $S_1(\oO,\oO)$.

In this section we study analogues of these cover selection properties for subsets of $\Kont^\kappa$.

The proofs in this section bear a strong resemblance to the known proofs for the case $\kappa=\omega$. The main new obstacle, which sometimes appears, is checking of the fact that there exists an appropriate subsequence of length $\kappa$ of a given sequence of covers.

\subsubsection{$\kappa$-$\gamma$-Sets}

A~family $\oU$ of open subsets  of a~topological space $X$ will be called a~\df{$\kappa$-cover} of $X$ if for any $A\in[X]^{<\kappa}$ there exists $U\in\oU$ such that $A\subseteq U$. It is a~\df{$\kappa$-$\gamma$-cover} if $\oU=\{U_\alpha\colon \alpha<\kappa\}$, and
$X\subseteq \bigcup_{\alpha<\kappa}\bigcap_{\alpha<\beta<\kappa} U_\beta$.
Notice that every subsequence of length $\kappa$ of a~$\kappa$-$\gamma$-cover is still a~$\kappa$-$\gamma$-cover.

The family of all $\kappa$-covers of $X$ will be denoted by $\Omega_\kappa (X)$, and the family of all $\kappa$-$\gamma$-covers will be denoted by $\Gamma_\kappa(X)$. The family of all open covers of size $\kappa$ of $X$, is denoted by $\oO_\kappa(X)$. The underlying set can be omitted in this notation if it is apparent from the context. We always assume that the covers which are considered in this paper are proper, i.e. the set itself is never an element of its cover.

 $X\subseteq \Kont^\kappa$ will be called a~\df{$\kappa$-$\gamma$-set} if for every open $\kappa$-cover $\oU$ of $X$ there exists a~sequence $\langle U_\alpha\rangle_{\alpha<\kappa}\in \oU^\kappa$ such that $\{U_\alpha\colon \alpha<\kappa\}$ is a~$\kappa$-$\gamma$-cover.

If $\oA,\oB$ are families of open covers of a~set $X$, we shall say that it satisfies the \df{$S^\kappa_1(\oA,\oB)$ property} if for every sequence $\langle \oU_\alpha\rangle_{\alpha<\kappa}\in\oA^\kappa$, there exists a~sequence $\langle U_\alpha\rangle_{\alpha<\kappa}$ such that $U_\alpha\in \oU_\alpha$, for all $\alpha<\kappa$, and $\{U_\alpha\colon \alpha<\kappa\}\in \oB$. 

Our aim is to prove that similarly to the case $\kappa=\omega$, $\kappa$-$\gamma$-sets can be characterized in terms of selection principles. First we need the following easy observation.

\begin{lemma}[$\omega$: \cite{lb:srl}, Theorem~8.95]\label{s1refinement}
Let $X$ be a~subset of a~$\kappa$-additive topological space, and $\oA,\oB$ be any families of open covers of cardinality $\kappa$ of $X$ such that:
\begin{enumerate}[(a)]
\item if in $\oB$ there exists a~refinement of an open cover $\oU$, then there also exists $\oU'\subseteq \oU$ with $\oU'\in \oB$,
\item if $\beta<\kappa$, and $\langle \oU_\alpha\rangle_{\alpha<\beta}\in \oA^\beta$, then there exists $\oU\in \oA$ such that $\oU$ is a~refinement of $\oU_\alpha$ for every $\alpha<\beta$,
\item if $\{U_\alpha\colon \alpha<\kappa\}\in \oB$,  $\oV_{\beta}=\{V_{\alpha,\beta}\colon \alpha<\gamma_{\beta}\}$ for $\beta<\kappa$ and $\langle \gamma_\beta\rangle_{\beta\in \kappa}\in \kappa^\kappa$ are such that $U_\beta\subseteq V_{\alpha,\beta}$ for all $\beta<\kappa, \alpha<\gamma_\beta$, then $\bigcup_{\beta<\kappa}\oV_\beta\in\oB$.
\end{enumerate}
Then $X$ satisfies $S^\kappa_1(\oA,\oB)$ if and only if for every $\langle \oU_\alpha\rangle_{\alpha<\kappa}\in \oA^\kappa$ such that $\oU_\beta$ is a~refinement of $\oU_\alpha$, for all $\alpha<\beta<\kappa$, there exists $\langle U_\alpha\rangle_{\alpha<\kappa}$ with $\{U_\alpha\colon \alpha<\kappa\}\in\oB$, and $U_\alpha\in\oU_\alpha$ for all $\alpha<\kappa$.
\end{lemma} 

Proof: 
Clear.\  \hfill $\square$

\begin{lemma}[$\omega$: \cite{lb:srl}]\label{lemgamma}
If $X\subseteq \Kont^\kappa$, then $\oA=\Omega_\kappa(X)$ and $\oB=\Gamma_\kappa(X)$ satisfy the premise of Lemma~\ref{s1refinement}.
\end{lemma} 

Proof: Recall that an intersection of less than $\kappa$ open sets in $\Kont^\kappa$ is still open. The rest of the proof is obvious.\hfill $\square$

\begin{theo}[$\omega$: \cite{lb:srl}, Theorem~8.96]\label{gamma1}
A~set $X\subseteq \Kont^\kappa$, with $|X|\geq\kappa$ is a~$\kappa$-$\gamma$-set if and only if it has property $S_1^\kappa(\Omega_\kappa,\Gamma_\kappa)$.
\end{theo}

Proof:   The proof is a relatively easy generalization of the proof in the case $\kappa=\omega$, which can be found in \cite[Theorem~8.96]{lb:srl}. In the last argument one has to apply Lemmas~\ref{s1refinement} and~\ref{lemgamma}.
 \hfill $\square$

\begin{cor}[$\omega$: \cite{lb:srl}, Corollary~8.97]
Every $\kappa$-$\gamma$-set satisfies $S_1^\kappa(\Gamma_\kappa,\Gamma_\kappa)$.
\end{cor}

Proof: Obviously, every $\kappa$-$\gamma$-cover is a~$\kappa$-cover. \hfill $\square$

Finally, we prove that every union of $\kappa$ many closed subsets of $\kappa$-$\gamma$-set is $\kappa$-$\gamma$-set as well.

We will also need the following fact, which in the case of $\kappa=\omega$ is noticed in \cite{lb:srl}[Exercise~8.21].

\begin{prop}[$\omega$: \cite{lb:srl}]\label{cl-gamma}
A~$\kappa$-union of closed subsets of a~$\kappa$-$\gamma$-set is a~$\kappa$-$\gamma$-set.
\end{prop}

Proof:  
Let $F=\bigcup_{\alpha<\kappa} F_\alpha$ with $F_\alpha\subseteq X$, where $X$ is a~$\kappa$-$\gamma$-set, and $F_\alpha$ are closed in $X$. Assume that for $\alpha<\beta<\kappa$, $F_\alpha\subseteq F_\beta$, and let $\oU$ be a~$\kappa$-cover of $F$. For any $\alpha<\kappa$, 
$\oU_\alpha=\{U\cup (X\setminus F_\alpha)\colon U\in\oU\}$
is a~$\kappa$-cover of $X$. Thus, by Theorem~\ref{gamma1}, there  exists a~sequence $\langle U_{\alpha}\rangle_{\beta<\kappa}$ such that $U_\alpha\in \oU_\alpha$, and $X\subseteq \bigcup_{\gamma<\kappa}\bigcap_{\gamma<\beta<\kappa} U_{\beta}$. Let $\langle V_{\alpha}\rangle_{\alpha,\beta<\kappa}\in \oU^\kappa$ be such that $U_{\alpha}=V_{\alpha}\cup (X\setminus F_\alpha)$. 

 Then $F\subseteq \bigcup_{\alpha<\kappa}\bigcap_{\alpha<\beta<\kappa} V_\beta,$ 
because if $x\in F$, then there exists $\alpha<\kappa$ such that $x\notin X\setminus F_\beta$ for all $\beta<\kappa$ with $\alpha<\beta$. Thus,
$x\in \bigcap_{\alpha<\beta<\kappa} V_\beta$.  \hfill $\square$

\subsubsection{$\kappa$-Hurewicz property}

A~cover $\oU$ of a~set $X$ is \df{essentially of size $\kappa$} if for every $\oV\in[\oU]^{<\kappa}$, $X\setminus \bigcup\oV\neq \zbp$.

We shall say that a~set $X$ satisfies the \df{$U_{<\kappa}^\kappa(\oA,\oB)$ principle} if for every sequence $\langle\oU_{\alpha}\rangle_{\alpha<\kappa}\in \oA^\kappa$ of covers essentially of size $\kappa$, there exists $\langle\oV_\alpha\rangle_{\alpha<\kappa}$ such that $\oV_\alpha\in[\oU_\alpha]^{<\kappa}$ for all $\alpha<\kappa$, and $\{\bigcup\oV_\alpha\colon \alpha<\kappa\}\in \oB$.

A~set $X$ has \df{$\kappa$-Hurewicz property} if it satisfies the $U_{<\kappa}^\kappa(\oO_\kappa,\Gamma_\kappa)$ principle.

\begin{prop}[$\omega$: \cite{lb:srl}, Theorem~8.98]\label{hurewicz1}
If $X\subseteq \Kont^\kappa$ satisfies $S^\kappa_1(\Gamma_\kappa,\Gamma_\kappa)$, then it has $\kappa$-Hurewicz property.
\end{prop}

Proof:   The proof is similar to the proof in the case $\kappa=\omega$, (see \cite[Theorem~8.98]{lb:srl}).
\hfill $\square$

\begin{samepage}
\begin{cor}[$\omega$: \cite{lb:srl}, Corollary~2.99]\label{hurewicz2}
If $X$ is a $\kappa$-$\gamma$-set, then it has $\kappa$-Hurewicz property.
\end{cor}
Proof: Clear. \hfill $\square$
\end{samepage}

Next we show that no Lusin set for $\kappa$ can have $\kappa$-Hurewicz property. We apply the following lemma.

\begin{lemma}[$\omega$: \cite{lb:srl}, Lemma~8.100]
If $A\subseteq \Kont^\kappa$ with an empty interior has $\kappa$-Hurewicz property, then $A$ is $\kappa$-meagre.
\end{lemma}

Proof: The proof is a relatively straightforward generalization of the proof in the case $\kappa=\omega$, which can be found in \cite[Lemma~8.100]{lb:srl}.
\hfill $\square$

\begin{cor}[$\omega$: \cite{lb:srl}, Theorem~8.101]\label{lusin-hur}
If $\kappa<\lambda\leq 2^\kappa$, and $L\subseteq \Kont^\kappa$ is a~$\lambda$-$\kappa$-Lusin set, then $L$ does not have $\kappa$-Hurewicz property.
\end{cor}
Proof: Clear. \hfill $\square$

\subsubsection{$\kappa$-Menger property}

A~set has \df{$\kappa$-Menger property} if it satisfies the $U_{<\kappa}^\kappa(\oO_\kappa,\oO_\kappa)$ principle.

Despite the fact that every Lusin set for $\kappa$ lacks $\kappa$-Hurewicz property (see Corollary~\ref{lusin-hur}), it has $\kappa$-Menger property.

\begin{prop}[$\omega$: \cite{lb:srl},  Theorem~8.101]\label{lus-men}
Let $L\subseteq \Kont^\kappa$ be a~Lusin set for $\kappa$. Then $L$ has $\kappa$-Menger property.
\end{prop}

Proof: The proof is a relatively straightforward generalization of the proof in the case $\kappa=\omega$, which can be found in \cite[Theorem~8.101]{lb:srl}.
\hfill $\square$ 

\subsubsection{$\kappa$-Rothberger property}

A~set has \df{$\kappa$-Rothberger property} if it satisfies the $S^\kappa_1(\oO_\kappa,\oO_\kappa)$ principle. Obviously, this property implies $\kappa$-Menger property.

\begin{prop}[$\omega$: \cite{lb:srl}, Theorem~8.102]\label{roth-sn}
If $A\subseteq \Kont^\kappa$ has $\kappa$-Rothberger property, then $A\in \SN_\kappa$.
\end{prop}

Proof: Clear. \hfill $\square$

\begin{cor}[$\omega$: \cite{lb:srl}, Corollary~8.103]
The generalized Cantor space $2^\kappa$ and the generalized Baire space $\kappa^\kappa$ do not have $\kappa$-Rothberger property.
\end{cor}
Proof: Clear. \hfill $\square$

Proposition~\ref{consn} can be formulated in a~stronger form.

\begin{prop}[$\omega$: \cite{lb:srl}, Theorem~8.105]\label{con-roth}
If $A\subseteq \Kont^\kappa$ is $\kappa^+$-concentrated on a~set $B\subseteq \Kont^\kappa$, $B\subseteq A$ with $|B|\leq \kappa$, then $A$ has $\kappa$-Rothberger property.
\end{prop}

Proof:  The proof is a generalization of the proof in the case $\kappa=\omega$ (see \cite[Theorem~8.105]{lb:srl}).
\hfill $\square$

This allows us to formulate a~stronger version of Proposition~\ref{lus-men}.

\begin{cor}[$\omega$: \cite{lb:srl}, Corollary~8.106]
Every Lusin set for $\kappa$ has $\kappa$-Rothberger property.
\end{cor}

Proof: By Proposition~\ref{lusin-con}, every Lusin set for $\kappa$ satisfies the premise of Proposition~\ref{con-roth}. \hfill $\square$

\begin{samepage}
\begin{lemma}[$\omega$: \cite{lb:srl}]\label{lemoo}
If $X\subseteq \Kont^\kappa$, then $\oA=\oO_\kappa(X)$ and $\oB=\oO_\kappa(X)$ satisfy the premise of Lemma~\ref{s1refinement}.
\end{lemma} 

Proof: Clear. \hfill $\square$
\end{samepage}

\begin{theo}[$\omega$: \cite{lb:srl}, Theorem~8.109]\label{rothgamma}
Every $\kappa$-$\gamma$-set of cardinality $\geq \kappa$ has $\kappa$-Rothberger property.
\end{theo}

Proof:  The proof is similar to the proof in the case $\kappa=\omega$, which can be found in \cite[Theorem~8.109]{lb:srl}. In the last argument one has to apply Lemmas~\ref{s1refinement} and~\ref{lemoo}.
\hfill $\square$

\begin{cor}[$\omega$: \cite{lb:srl}, Corollary~8.110]
Every $\kappa$-$\gamma$-set is $\kappa$-strongly null.
\end{cor}

Proof: Follows by Corollary~\ref{roth-sn}. \hfill $\square$

\begin{cor}
The generalized Cantor space $2^\kappa$ and the generalized Baire space $\kappa^\kappa$ are not $\kappa$-$\gamma$-sets.
\end{cor}

Proof: Clear. \hfill $\square$

Thus, no $\kappa$-perfect subset of $2^\kappa$ is a $\kappa$-$\gamma$-set. Nevertheless, the following question remains unanswered.

\begin{prob}
Is there a closed subset of $2^\kappa$ which is a $\kappa$-$\gamma$-set?
\end{prob}

We finish this part by proving a~lemma, which becomes useful in the next section.

\begin{lemma}[$\omega$: \cite{antw:rpsssas}, Theorem~2.1 (1)]\label{lem-gamma}
Assume that $\kappa$ is a~weakly inaccessible  cardinal. Let $A\subseteq 2^\kappa$ be a~$\kappa$-$\gamma$ set, which is not closed. Then there exists $B\in [\kappa]^\kappa$ such that for all $C\in [B]^\kappa$, $\chi_C\notin A$.
\end{lemma}

Proof: Let $A\subseteq 2^\kappa$ be a~$\kappa$-$\gamma$ set, and let $b\colon\bigcup_{\alpha<\kappa}\{\alpha\}\times \alpha\to \kappa$ be a~bijection. Notice that $2^\kappa\setminus A$ is not an open set, as $A$ is not closed. Therefore, there exists $y\in 2^\kappa\setminus A$ such that $A\cap[y\obc\alpha]\neq\zbp$, for any $\alpha<\kappa$. Choose inductively a~sequence $\langle x_\alpha\rangle_{\alpha<\kappa}\in A^\kappa$ such that  for $\alpha,\beta<\kappa$, if $x_\alpha=x_\beta$, then $\alpha=\beta$, and such that for every $\gamma<\kappa$ there exists $\alpha<\kappa$ such that $y\obc \gamma = x_\alpha\obc \gamma$. To achieve this take any $x_0\in A$, and for $\alpha<\kappa$, let 
$\xi = \bigcup_{\beta<\alpha}\bigcup\{\gamma<\kappa\colon y\obc\gamma = x_\beta\obc\gamma\}$.
Let $x_\alpha\in A\cap [y\obc\xi+1]$.

If $I\subseteq \kappa$ and $s\in 2^I$, assume that $[s]$ denote $\{x\in 2^\kappa\colon x\obc I=s\}$. For $\alpha<\kappa$, let
\[\oU_{\alpha}=\left\{\left(\bigcup_{s\in S} [s]\cap A\setminus \bigcup_{\alpha\leq \beta<\kappa} \{x_{\beta}\}\right)\colon S\in \left[2^{b[\{\alpha\}\times \alpha]}\right]^{<|\alpha|}\right\},\]
and define $\oU=\bigcup_{\alpha<\kappa}\oU_\alpha$. Notice that $\oU$ is a~$\kappa$-cover of $A$, because $\kappa$ is weakly inaccessible. Therefore, there is $\langle U_\alpha\rangle_{\alpha<\kappa}\in \oU^\kappa$ such that
$A\subseteq \bigcup_{\alpha<\kappa}\bigcap_{\alpha<\beta<\kappa} U_\beta$.

But since 
$x_\alpha\notin \bigcup_{\beta<\alpha}\bigcup\oU_\beta,$
 for all $\alpha< \kappa$, we get that for any $\alpha<\kappa$, there is $\xi<\kappa$ such that for all $\xi<\beta<\kappa$, there is $\alpha<\gamma<\kappa$, so that $U_\beta\in\oU_\gamma$. Therefore, we can choose inductively increasing sequences $\langle\xi_\alpha\rangle_{\alpha<\kappa}\in \kappa^\kappa$ and $\langle\delta_\alpha\rangle_{\alpha<\kappa}\in \kappa^\kappa$ such that $U_{\xi_\alpha}\in \oU_{\delta_\alpha}$, for any $\alpha<\kappa$. 
 
Fix $\alpha<\kappa$, and let $S_\alpha\in\left[2^{b[\{\delta_\alpha\}\times \delta_\alpha]}\right]^{<|\delta_\alpha|}$ be such that 
$U_{\xi_\alpha}=\bigcup_{s\in S_\alpha}[s]\cap A\setminus  \bigcup_{\delta_\alpha\leq \beta<\kappa} \{x_{\beta}\}$.
There exists $\eta_\alpha<\delta_\alpha$ such that $\chi_{a}\notin S_\alpha$ for any $a\supseteq \{b(\delta_\alpha,\eta_\alpha)\}$. Put
$B=\{b(\delta_\alpha,\eta_\alpha)\colon \alpha<\kappa\}$.
Then, for all $C\in [B]^\kappa$, $\chi_C\notin A$. Indeed, if $C\in [B]^\kappa$, then for every $\alpha<\kappa$, there is $\alpha<\beta<\kappa$ such that for
$a=C\cap b[\{\delta_\beta\}\times \delta_\beta\}]=\{b(\delta_\beta,\eta_\beta)\}$
we have that $\chi_a\notin S_\beta$. 
For such $\beta$, $\chi_C\notin U_{\xi_\beta}$, therefore for all $\alpha<\kappa$,
$\chi_C\notin \bigcap_{\alpha<\beta<\kappa}U_{\xi_\beta}$,
and hence
$\chi_C\notin \bigcup_{\alpha<\kappa}\bigcap_{\alpha<\beta<\kappa} U_\beta \supseteq A$.

\ \hfill $\square$

\section{Generalization of other notions of smallness in $2^\kappa$ and $\kappa^\kappa$}\label{chother}

In this section, we present generalizations of some less common notions of small sets.

Some of the results presented here have their counterparts in the standard case of $2^\omega$ (or $\omega_1^{\omega_1}$), and if so, we give a~reference in the form ($\omega$: [n]) (or ($\omega_1$: [n])).  

In particular, we consider Ramseyan properties of some special subsets of $\Kont^\kappa$ in terms of the tree-like forcings that generalize forcings which appeared in the first part of this paper. These notions and their natural variants were investigated in  papers\cite{sfykvk:rpgr} and \cite{gl:gssm}. In \cite{sfykvk:rpgr}, the authors study regularity properties of Borel subsets of $\Kont^\kappa$, and they point out that the behaviour of $\Sigma_1^1$ subsets in relation to the generalized tree-like forcing is different than this for $\kappa=\omega$. In this part, we are mainly interested is relations between these ideals and the classes of sets discussed in the previous section.

\subsection{$X$-small sets}
In this section, we present some parallels to the results from \cite[Chapter 4]{ah:nsgbs}.

If $X\subseteq \kappa$, then a~set $A\subseteq \Kont^\kappa$ will be called \df{$X$-small} if there exists $\langle a_\alpha\rangle_{\alpha\in X}\in (\Kont^\kappa)^X$ such that
$A\subseteq \bigcup_{\alpha\in X}[a_\alpha\obc \alpha]$.
Notice that $A$ is $\SN_\kappa$ if it is $X$-small for any $X\in [\kappa]^\kappa$.

Consider the following order on $[\kappa]^\kappa$. For $X,Y\subseteq [\kappa]^\kappa$, let $X<Y$ (respectively, $X\leq Y$) if and only if there exists a~bijection $F\colon X\to Y$ such that for all $\alpha\in X$, $\alpha< F(\alpha)$ (respectively, $\alpha\leq F(\alpha)$). Let $X+1=\{\alpha+1\colon \alpha\in X\}$. Notice that if $X<Y$, then $X+1\leq Y$.

Let $X,Y\in [\kappa]^\kappa$ be such that $X<Y$. Then, the family of $Y$-small sets is a~proper subfamily of $X$-small sets (see \cite{ah:nsgbs}, Proposition~4.2).

Let $\lambda<\kappa$. We say that a~set $A\subseteq \Kont^\kappa$ is \df{$\lambda$-$X$-small} for $X\subseteq \kappa$ if there exists $\langle a_{\alpha,\beta}\rangle_{\alpha\in X, \beta<\lambda}\in (\Kont^\kappa)^{X\times \lambda}$ such that $A\subseteq \bigcup_{\alpha\in X}\bigcup_{\beta<\lambda}[a_{\alpha,\beta}\obc \alpha]$.
We say that $A\subseteq \Kont^\kappa$ is \df{$\oX$-null} for $\oX\subseteq [\kappa]^{\leq\kappa}$ if for all $X\in \oX$, $A$ is $X$-small, and \df{$\lambda$-$\oX$-null} if for all $X\in \oX$, $A$ is $\lambda$-$X$-small.

The notion of $\lambda$-$\oX$-null sets for $\oX\subseteq [\kappa]^\lambda$ does not depend precisely on $\oX$. Indeed, we get the following proposition.

\begin{prop}\label{null}
Let $\lambda<\kappa$. A~set $A\subseteq \Kont^\kappa$ is $\lambda$-$\{\{\alpha\}\colon \alpha<\kappa\}$-null in $\Kont^\kappa$ if and only if it is $[\kappa]^\lambda$-null.
\end{prop}

Proof: Let $\lambda<\kappa$, and assume that $A$ is $\lambda$-$\{\{\alpha\}\colon \alpha<\kappa\}$-null. Let $X=\{\xi_\beta\colon\beta<\lambda\}\in [\kappa]^\lambda$ and $\alpha=\bigcup X$. Obviously, $\alpha<\kappa$. Therefore, there exists a~sequence $\langle a_\beta\rangle_{\beta<\lambda}$ such that
$A\subseteq \bigcup_{\beta<\lambda} [a_{\beta}\obc \alpha]\subseteq \bigcup_{\beta<\lambda} [a_{\beta}\obc \xi_{\beta}]$,
so $A$ is $X$-small. 

On the other hand, assume that $A$ is $[\kappa]^\lambda$-null and $\alpha<\kappa$. Let $X=\{\alpha+\beta\colon \beta<\lambda\}\in [\kappa]^\lambda$. There exists a~sequence $\langle a_\beta\rangle_{\beta<\lambda}$ such that
\[A\subseteq \bigcup_{\beta<\lambda} [a_{\beta}\obc \alpha+\beta]\subseteq \bigcup_{\beta<\lambda} [a_{\beta}\obc \alpha],\]
so $A$ is $\lambda$-$\{\alpha\}$-small. \ \hfill $\square$ 

A~set $A\subseteq \Kont^\kappa$ will be called \df{small in $\Kont^\kappa$} if 
there exists $\lambda<\kappa$ such that 
$A$ is $\lambda$-$\{\{\alpha\}\colon \alpha<\kappa\}$-null. Obviously, every $A\subseteq \Kont^\kappa$ with $|A|<\kappa$ is small in $\Kont^\kappa$.

Notice that every small set in $\Kont^\kappa$ is $\kappa$-strongly null.

\begin{prop}[$\omega_1$: \cite{ah:nsgbs}, Lemma~4.3 ii)]\label{smallsn}
Assume that a set $A\subseteq \Kont^\kappa$ is small in $\Kont^\kappa$. Then $A\in\SN_\kappa$.
\end{prop}

Proof:
 Let $\lambda<\kappa$ be such that $A$ is $\lambda$-$\{\{\alpha\}\colon \alpha<\kappa\}$-null. Therefore, by Proposition~\ref{null}, $A$ is $[\kappa]^\lambda$-null. Let $X=\{\xi_{\alpha}\colon \alpha<\kappa\}\in [\kappa]^\kappa$. There exists a~sequence $\langle a_{\alpha}\rangle_{\alpha<\lambda}\in (\Kont^\kappa)^\lambda$ such that $A\subseteq \bigcup_{\alpha<\lambda}[a_{\alpha}\obc \xi_\alpha]$. For $\lambda\leq \alpha<\kappa$ set $a_\alpha=\bzero$. We get that $A\subseteq \bigcup_{\alpha<\kappa}[a_{\alpha}\obc \xi_\alpha]$.
\ \hfill $\square$

\begin{prop}[$\omega_1$: \cite{ah:nsgbs}, Proposition~4.4]
A~set $A\subseteq \Kont^\kappa$ is $\SN_\kappa$ if and only if there exists $\lambda<\kappa$ such that $A$ is $\lambda$-$[\kappa]^\kappa$-null.
\end{prop}

Proof: If $A\subseteq \Kont^\kappa$ is $\SN_\kappa$, it is  obviously $\lambda$-$[\kappa]^\kappa$-null for all $\lambda<\kappa$.  Assume that $\lambda<\kappa$, and $A\subseteq \Kont^\kappa$ is $\lambda$-$[\kappa]^\kappa$-null. Let $X=\{\xi_\alpha\colon \alpha<\kappa\}\in [\kappa]^\kappa$. Let $b\colon \lambda\times\kappa\to\kappa$ be a~bijection, and for all $\alpha<\kappa$, let $X_\alpha=\{\xi_{b(\beta,\alpha)}\colon\beta<\lambda\}\in [\kappa]^\lambda$. Put $\delta_\alpha=\bigcup X_\alpha$, for $\alpha<\kappa$. Finally, let $Y=\{\delta_\alpha\colon \alpha<\kappa\}\in [\kappa]^\kappa$. We can find $\langle x_{\alpha,\beta}\rangle_{\alpha<\kappa,\beta<\lambda}\in (\Kont^\kappa)^{\kappa\times\lambda}$ such that
$A\subseteq \bigcup_{\alpha<\kappa}\bigcup_{\beta<\lambda}[x_{\alpha,\beta}\obc\delta_\alpha]$.
For $\alpha<\kappa$, let $z_{\alpha}=x_{b^{-1}(\alpha)}$. Then
$A\subseteq \bigcup_{\alpha<\kappa}[z_{\alpha}\obc\delta_{\pi_2(b^{-1}(\alpha))}]\subseteq \bigcup_{\alpha<\kappa}[z_{\alpha}\obc\xi_{\alpha}]$,
where $\pi_2$ is the projection on $\kappa$.
\hfill $\square$

\begin{prop}
Let $X\subseteq\kappa$ be such that $0\notin X$ and $X\cap\Lim=\zbp$. If $A\subseteq 2^\kappa$ is $X$-small, then $|2^\kappa\setminus A|=2^\kappa$.
\end{prop}

Proof: Let $\langle x_\alpha\rangle_{\alpha\in X}\in \left(2^\kappa\right)^X$ be such that $A\subseteq\bigcup_{\alpha\in X}[x_\alpha\obc \alpha]$. Consider the~set
\[B=\left\{x\in 2^\kappa\colon \forall_{\alpha<\kappa}\left(\alpha+1\in X\Rightarrow x(\alpha)=x_{\alpha+1}(\alpha)+1\right)\right\}.\]
Then for all $\alpha\in X$, $B\cap [x_\alpha\obc\alpha]=\zbp$. Thus, $B\cap A=\zbp$. Furthermore, $B$ contains a set homeomorphic to $2^\kappa$, so $|2^\kappa\setminus A|=2^\kappa$. \hfill $\square$

Next we study a~connection between the diamond principle for $\kappa$ (see section~\ref{intro-gen}) and the notion of $C$-smallness for closed unbounded or stationary sets $C\subseteq \kappa$.

For $E\subseteq \kappa$, $A\subseteq 2^\kappa$, $\oI\subseteq [\kappa]^{\leq\kappa}$, let \df{$\diamondsuit_\kappa(E,A,\oI)$} denote the following principle: there exists a~sequence $\langle s_\alpha\rangle_{\alpha<\kappa}\in \left(2^{<\kappa}\right)^\kappa$ such that for all $x\in A$, $\{\alpha\in E\colon x\obc \alpha=s_{\alpha}\}\notin \oI$.

Notice the following easy observation.

\begin{prop}
If $A\subseteq 2^\kappa$, and $E\subseteq \kappa$, then $\diamondsuit_\kappa(E,A,\{\zbp\})$ if and only if $A$ is $E$-small.
\end{prop}

Proof: Indeed, $\diamondsuit_\kappa(E,A,\{\zbp\})$ if and only if 
for all $x\in A$, $\{\alpha\in E\colon x\obc \alpha=s_{\alpha}\}\neq\zbp$. \hfill $\square$

Let \df{$\NS_\kappa$} be the family of all non-stationary sets in $\kappa$. 

\begin{prop}[$\omega_1$: \cite{ah:nsgbs}, Lemma~5.3]\label{diamond}
Let $E\subseteq \kappa$. The principle $\diamondsuit_\kappa(E,2^\kappa,\NS_\kappa)$ holds if and only if $\diamondsuit_\kappa(E)$ holds.
\end{prop}

Proof: The proof is a relatively straightforward generalization which makes use of Proposition~\ref{diamond} of the proof in the case $\kappa=\omega_1$ (see \cite{ah:nsgbs}[Lemma~5.3]).
 \hfill $\square$

\begin{prop}[$\omega_1$: \cite{ah:nsgbs}, Proposition~4.6 iii)]\label{clubsmall}
Assume $\diamondsuit_\kappa$. If $C$ is a~closed unbounded set in $\kappa$, then $2^\kappa$ is $C$-small.
\end{prop}  

Proof:  The proof is a relatively straightforward generalization which makes use of Proposition~\ref{diamond} of the proof in the case $\kappa=\omega_1$,  (see \cite{ah:nsgbs}[Proposition~4.6 iii)]).
 \hfill $\square$
 
\begin{prop}[$\omega_1$: \cite{ah:nsgbs}, Proposition~4.6 iv)]\label{diamondEsmall}
Let $E\subseteq\kappa$, and assume $\diamondsuit_\kappa(E)$. Then $2^\kappa$ is $E$-small.
\end{prop}

Proof:  The proof is a relatively straightforward generalization which applies Proposition~\ref{diamond} of the proof in the case $\kappa=\omega_1$ (see \cite{ah:nsgbs}[Proposition~4.6 iv)]).
 \hfill $\square$

\begin{cor}[$\omega_1$: \cite{ah:nsgbs}, Proposition~4.6 iv)]\label{vlsmall}
Assume $V=L$. Then $2^\kappa$ is $X$-small for every stationary set $X\subseteq\kappa$.
\end{cor}

Proof: Recall that $V=L$ implies $\diamondsuit_\kappa(X)$ for every stationary set $X\subseteq \kappa$ (see \cite[Exercise VI.14]{kk:stiip}). Therefore, by Proposition~\ref{diamondEsmall}, $2^\kappa$ is small for every stationary $X\subseteq \kappa$. \hfill $\square$

The whole space $2^\kappa$ can be presented as a~union of a~$\kappa$-meagre set, and a~$\oX$-null set for $\oX\in \left[[\kappa]^\kappa\right]^\kappa$.

\begin{prop}[$\omega_1$: \cite{ah:nsgbs}, Proposition~4.9]
Let $\oX\in \left[[\kappa]^\kappa\right]^\kappa$. There exist $A,B\subseteq 2^\kappa$ such that $A$ is $\oX$-null, and $B$ is $\kappa$-meagre, and $A\cup B=2^\kappa$.
\end{prop}

Proof:  The proof is similar to the proof in the case $\kappa=\omega_1$, which can be found in \linebreak \cite{ah:nsgbs}[Proposition~4.9].
\hfill $\square$

On the other hand, we have the following.

\begin{prop}[$\omega_1$: \cite{ah:nsgbs}, Proposition~4.11]\label{smnd}
Every small set in $\Kont^\kappa$ is nowhere dense.
\end{prop}

Proof: Let $\lambda<\kappa$ be such that $A\subseteq 2^\kappa$ is $\lambda-\{\{\alpha\}\colon \alpha<\kappa\}$-null. Let $s\in 2^\beta$ with $\beta<\kappa$, and let $\xi=\beta+\lambda$. There exists $\langle x_\alpha\rangle_{\alpha<\lambda}\in (\Kont^\kappa)^\lambda$ such that $A\subseteq \bigcup_{\alpha<\lambda}[x_\alpha\obc\xi]$.
But $|\{x\obc\xi\colon x\in [s]\}|=2^\lambda$, thus there exists $t\in 2^\xi$ such that $s\subseteq t$, and $[t]\cap A=\zbp$. \hfill $\square$

However, not every nowhere dense set in $\Kont^\kappa$ is small in $\Kont^\kappa$.

\begin{prop}[$\omega_1$: \cite{ah:nsgbs}, Proposition~4.12]\label{ndnsm}
There exists a~nowhere dense set $A\subseteq \Kont^\kappa$ which is not $\kappa$-strongly null.
\end{prop}

Proof: Let $\langle \xi_\alpha\rangle\in \kappa^\kappa$ be an increasing sequence of limit ordinals. Let 
$A=\left\{x\in \Kont^\kappa\colon \forall_{\alpha<\kappa}x(\xi_\alpha)=0\right\}$.
Obviously, $A$ is nowhere dense. Assume that $A\in \SN_\kappa$. Then there exists $\langle x_\alpha\rangle_{\alpha<\kappa}\in (\Kont^\kappa)^\kappa$ such that 
$A\subseteq \bigcup_{\alpha<\kappa} [x_\alpha\obc \xi_\alpha+2]$. Let $x\in 2^\kappa$ be such that $x(\xi_\alpha+1)=x_\alpha(\xi_\alpha+1)+1$ for all $\alpha<\kappa$, and $x(\beta)=0$ for $\beta\notin\{\xi_\alpha+1\colon \alpha\in \kappa\}$. Then $x\in A$, but $x\notin \bigcup_{\alpha<\kappa} [x_\alpha\obc \xi_\alpha+2]$, which is a~contradiction. \hfill $\square$ 

\subsection{$\kappa$-Meagre additive sets}

In this section, we present some generalizations of results related to meagre additive sets. We start by noticing that  the combinatorial characterization of meagre sets (see \cite[Theorem 2.2.4]{tbhj:stsrl}) can be easily generalized. This generalization can be found in  \cite[Observation~5.1]{ss:pnii}. We present it here along with the proof for the sake of completeness.

\begin{prop}[\cite{ss:pnii}, Observation~5.1, $\omega$: \cite{tbhj:stsrl}, Theorem~2.2.4]\label{char-meagre}
Assume that $\kappa$ is strongly inaccessible, and $A\subseteq 2^\kappa$ is a~$\kappa$-meagre set. Then there exist an increasing sequence $\langle \xi_\alpha\rangle_{\alpha<\kappa}\in \kappa^\kappa$ and $y\in 2^\kappa$ such that $A\subseteq \left\{x\in 2^\kappa\colon \exists_{\beta<\kappa}\forall_{\beta<\gamma<\kappa} \exists_{\xi_{\gamma}\leq\xi<\xi_{\gamma+1}} x(\xi)\neq y(\xi)\right\}$.
\end{prop}

Proof: The proof is similar to the proof of \cite{tbhj:stsrl}[Theorem~2.2.4]. We use the fact that $\kappa$ is strongly inaccessible to find below $\kappa$ a bound  of $2^\xi$ many ordinals in the inductive step, where $\xi<\kappa$.
\ \hfill $\square$

Recall that a~set $A$ is called \df{meagre-additive} ($A\in \MM^{*}$) if for any meagre set $X$, $A+X$ is  meagre (see e.g. \cite{tw:manascs} and \cite{tbhj:stsrl}). The following \df{characterization of meagre-additive sets} is well-known. A set $X\in \MM^*$ (\cite{tbhj:stsrl}[Theorem 2.7.17]) if and only if for every increasing $f\in \omega^\omega$, there exists increasing $g\in \omega^\omega$ and $y\in 2^\omega$ such that for all $x\in X$, there exists $m\in \omega$, so that for each $n>m$, there is $k_n\in\omega$ with $g(n)\leq f(k_n)<f(k_n+1)\leq g(n+1)$ and such that 
$x\obc[f(k_n), f(k_n+1))=y\obc[f(k_n), f(k_n+1))$.

A~set $A\subseteq 2^\kappa$ will be called \df{$\kappa$-meagre additive} if for any $\kappa$-meagre set $F$, $A+F$ is $\kappa$-meagre. The family of all $\kappa$-meagre additive sets is denoted by $\MM_{\kappa}^*$. 

By Proposition~\ref{gms1}, we immediately get the following corollary.

\begin{cor}
Every $\kappa$-meagre additive set is $\kappa$-strongly null.
\end{cor}
Proof: Let $A$ be $\kappa$-meagre additive set. Since for every nowhere dense $F$, $A+F$ is $kappa$-meagre, we get that $A+F\neq 2^\kappa$, so there exists $x\in 2^\kappa$ such  that $(x+A)\cap F=\zbp$. Thus by Proposition~\ref{gms1}, $A$ is $\kappa$-strongly null.\ \hfill $\square$

The following theorem is a generalization of the characterization of meagre-additive sets (\cite[Theorem 2.7.17]{tbhj:stsrl}).

\begin{prop}[$\omega$: \cite{tbhj:stsrl}, Theorem 2.7.17]\label{madditive}
Assume that $\kappa$ is strongly inaccessible, and $X\subseteq 2^\kappa$. Then $X\in \MM_\kappa^*$ if and only if for every increasing sequence $\langle\xi_\alpha\rangle_{\alpha<\kappa}\in \kappa^\kappa$ there exist a~sequence $\langle\eta_\alpha\rangle_{\alpha<\kappa}\in \kappa^\kappa$ and $z\in 2^\kappa$ such that
\[X\subseteq \left\{ x\in 2^\kappa \colon \exists_{\alpha<\kappa}\forall_{\alpha<\beta<\kappa}\exists_{\gamma<\kappa} \left(\eta_\beta\leq \xi_\gamma <\xi_{\gamma+1}\leq \eta_{\beta+1}\land \forall_{\xi_\gamma\leq \delta< \xi_{\gamma+1}} x(\delta)=z(\delta)\right)\right\}.\]
\end{prop} 

Proof:  The proof is a  straightforward generalization  of the proof in the case $\kappa=\omega$, which can be found in \cite{tbhj:stsrl}[Theorem 2.7.17]. The assumption that $\kappa$ is strongly inaccessible is used here to apply Proposition~\ref{char-meagre}. 
\ \hfill $\square$

Notice that this implies that under the same assumption every $\kappa$-meagre additive set is $\PkM_\kappa$.

\begin{prop}[$\omega$: \cite{mkantw:ssrtfn}, Proposition~3.4]\label{mapkmk}
Assume that $\kappa$ is a~strongly inaccessible cardinal. Then every $\kappa$-meagre additive set is $\PkM_\kappa$.
\end{prop}

Proof: Let $A\in\MM_\kappa^*$, and let $P\subseteq 2^\kappa$ be a~$\kappa$-perfect set. By induction we construct a~sequence $\langle\xi_\alpha\rangle_{\alpha<\kappa}\in \kappa^\kappa$ such that $\xi_0=0$, and for $\alpha<\kappa$,
$\xi_{\alpha+1}=\bigcup_{t\in T_P\cap 2^{\xi_\alpha}}\min\left\{\len(s)\colon t\subseteq s\in \Spl(T_P)\right\}+1$.
Finally, for limit $\alpha<\kappa$, let $\xi_\alpha=\bigcup_{\beta<\alpha} \xi_\beta$. 

By Proposition~\ref{madditive}, we can find  a~sequence $\langle\eta_\alpha\rangle_{\alpha<\kappa}\in \kappa^\kappa$ and $z\in 2^\kappa$ such that
\[A\subseteq \bigcup_{\alpha<\kappa}\left\{ x\in 2^\kappa \colon \forall_{\alpha<\beta<\kappa}\exists_{\gamma<\kappa} \left(\eta_\beta\leq \xi_\gamma <\xi_{\gamma+1}\leq \eta_{\beta+1}\land \forall_{\xi_\gamma\leq \delta< \xi_{\gamma+1}} x(\delta)=z(\delta)\right)\right\}.\]

Let $\alpha<\kappa$, and let $s\in T_P$. Fix $s'\in T_P$ such that $s\subseteq s'$, and for some $\beta>\alpha$, $\len(s')=\eta_\beta$. 
Define $\gamma_0=\min\left\{\gamma<\kappa\colon\eta_\beta\leq \xi_\gamma <\xi_{\gamma+1}\leq \eta_{\beta+1}\right\}$
and
$\gamma_1= \bigcup\left\{\gamma<\kappa\colon  \eta_\beta\leq \xi_\gamma <\xi_{\gamma+1}\leq \eta_{\beta+1}\right\}+1$.

We construct  a~sequence $\langle t_{\delta}\rangle_{\gamma_0\leq\delta\leq\gamma_1}$ inductively, so that for all $\gamma_0\leq\delta\leq\delta'\leq\gamma_1$, $t_{\delta}\in T_P\cap 2^{\xi_\delta}$, $t_\delta\subseteq t_{\delta'}$, and $\exists_{\xi_\delta\leq \xi< \xi_{\delta+1}} t_{\delta+1}(\xi)\neq z(\xi)$.
Indeed, let $t_{\gamma_0}\in T_P$ be such that $s\subseteq t_{\gamma_0}$, and $\len(t_{\gamma_0})=\xi_{\gamma_0}$. Given $t_\delta$, by definition of $\langle \xi_\alpha\rangle_{\alpha<\kappa}$, one can find $t_{\delta+1}\supseteq t_\delta$ such that $\exists_{\xi_\delta\leq \xi< \xi_{\delta+1}} t_{\delta+1}(\xi)\neq z(\xi)$, as $|\{t\in T_P\cap 2^{\xi_{\delta+1}}\colon t\supseteq t_{\delta}\}|\geq 2$. For limit $\delta<\kappa$, choose any $t_\delta\supseteq \bigcup_{\gamma_0\leq\xi<\delta} t_\xi$ such that $\len(t_\delta)=\xi_\delta$.

Then,
 \[[t_{\gamma_1}]\cap P\cap \left\{ x\in 2^\kappa \colon \forall_{\alpha<\beta<\kappa}\exists_{\gamma<\kappa} \left(\eta_\beta\leq \xi_\gamma <\xi_{\gamma+1}\leq \eta_{\beta+1}\land \forall_{\xi_\gamma\leq \delta< \xi_{\gamma+1}} x(\delta)=z(\delta)\right)\right\}\]
is empty, and hence $A$ is $\kappa$-meagre in $P$. \hfill $\square$

We do not know whether the above result can be obtained for other cardinals, which satisfy weaker assumptions.
 
\begin{prob}
Assume that $\kappa$ is a cardinal which is not strongly inaccessible. Is every $\kappa$-meagre additive set in $\PkM_\kappa$?
\end{prob}

\subsection{$\kappa$-Ramsey null sets}\label{s-rn}

A~set $A$ is \df{Ramsey null} (see \cite{sp:crs}) if for any $n\in\omega$, $s\in 2^n$ and $S\in [\omega\setminus n]^\omega$, there exists $S'\in [S]^\omega$ such that $[s,S']\cap A=\zbp$, where for $s\in 2^n$, $n\in\omega$ and $S\in [\omega\setminus n]^\omega$, 
\[[s,S]=\{x\in 2^\omega\colon s^{-1}[\{1\}]\subseteq x^{-1}[\{1\}]\subseteq s^{-1}[\{1\}]\cup S\land |x^{-1}[\{1\}]\cap S|=\omega\}.\]

In this section, we generalize some results presented in \cite{antw:rpsssas}.

For $\alpha<\kappa$, $s\in 2^\alpha$ and $S\in [\kappa\setminus \alpha]^\kappa$, let
\[[s,S]=\{x\in 2^\kappa\colon s^{-1}[\{1\}]\subseteq x^{-1}[\{1\}]\subseteq s^{-1}[\{1\}]\cup S \land |x^{-1}[\{1\}]\cap S|=\kappa\}.\]
A~set $A\subseteq 2^\kappa$ will be called \df{$\kappa$-Ramsey null ($\kappa-CR_0$)} if for any $\alpha<\kappa$, $s\in 2^\alpha$ and $S\in [\kappa\setminus \alpha]^\kappa$, there exists $S'\in [S]^\kappa$ such that $[s,S']\cap A=\zbp$.

It is a~well-known fact that the ideal of Ramsey null subsets of $2^\omega$ is a~$\sigma$-ideal (see e.g. \cite{lh:cst}). We do not know whether the analogue holds for $\kappa$-Ramsey null sets.

\begin{prob}\label{qrams}
Is the ideal of $\kappa$-Ramsey null subsets of $2^\kappa$ $\kappa^+$-complete?
\end{prob}

\begin{theo}[$\omega$: \cite{antw:rpsssas}, Theorem~2.1]\label{gammarams}
Assume that $\kappa$ is a~weakly inaccessible cardinal. Then every $\kappa$-$\gamma$-set which is not closed in $2^\kappa$ is $\kappa$-Ramsey null.
\end{theo}

Proof: The proof is based on the proof of \cite[Theorem~2.1]{antw:rpsssas}. Namely, let $A\subseteq 2^\kappa$ be a~$\kappa$-$\gamma$-set, let $\delta<\kappa$, $s\in 2^\delta$, and $S=\{\xi_\alpha\colon \alpha<\kappa\}\in [\kappa\setminus \delta]^\kappa$. Define 
$E=\{x\in 2^\kappa\colon s^{-1}[\{1\}]\subseteq x^{-1}[\{1\}]\subseteq s^{-1}[\{1\}]\cup S\}=s_0+S_0$,
where $s_0=s\cup\{\langle\beta,0\rangle\colon \beta\in\kappa\setminus\delta\}$, and $S_0=\{f\cup\{\langle\beta,0\rangle\colon \beta\notin S\}\colon f\in 2^S\}$. Notice that $S_0$ is a~closed set in $2^\kappa$, and so is $E$.  Moreover, $\varphi\colon 2^\kappa\to E$ given by the following expression
$\varphi(x)=s_0+\chi_{\{\xi_\alpha\colon x(\alpha)=1\land\alpha<\kappa\}}$
is a~homeomorphism.

By Proposition~\ref{cl-gamma}, $E\cap A$ is a $\kappa$-$\gamma$ set, and therefore so is $\varphi^{-1}[E\cap A]$. By Lemma~\ref{lem-gamma}, there exists $B\in [\kappa]^\kappa$ such that for all $C\in [B]^\kappa$, $\chi_C\notin \varphi^{-1}[E\cap A]$, which means that $\varphi(\chi_C)\notin A$. Let $S'=\{\xi_\alpha\colon \alpha\in B\}$. Then  $S'\in [S]^\kappa$, and $[s,S']=\{\varphi(\chi_C)\colon C\in [B]^\kappa\}$. Thus, $[s,S']\cap A=\zbp$. \hfill $\square$

\begin{lemma}[$\omega$: \cite{antw:rpsssas}, Lemma~4.1]\label{lem-mm}
If $A,B\subseteq 2^\kappa$, then
$2^\kappa\setminus\left(A+(2^\kappa\setminus B)\right)=\{x\in 2^\kappa\colon x+A\subseteq B\}$.
\end{lemma}

Proof: The proof of \cite{antw:rpsssas}[Lemma~4.1] is valid for any vector space over $\Z_2$.\hfill $\square$

\begin{prop}[$\omega$: \cite{antw:rpsssas}, Lemma~4.2]\label{mmrn}
Assume that $\kappa$ is strongly inaccessible, and $A\subseteq 2^\kappa$ is a~$\kappa$-meagre set. Then there exists a~$\kappa$-meagre set $B\subseteq 2^\kappa$ such that $A+(2^\kappa\setminus B)$ is $\kappa$-Ramsey null.
\end{prop}

Proof: The proof is a relatively straightforward generalization  of the proof in the case $\kappa=\omega$, which can be found in \cite{antw:rpsssas}[Lemma~4.2]. The assumption that $\kappa$ is strongly inaccessible is used to apply Proposition~\ref{char-meagre}. 
\hfill $\square$

We get the following theorem.

\begin{theo}[$\omega$: \cite{antw:rpsssas}, Theorem~4.3]\label{manrams}
Assume that $\kappa$ is strongly inaccessible, $\cov(\kappa-CR_0)= 2^\kappa$, and $\add(\MM_\kappa)=2^\kappa$. Then there exists a~$\kappa$-meagre additive set which is not $\kappa$-Ramsey null.
\end{theo}

Proof: Let $\{F_\alpha\colon\alpha<2^\kappa\}$ be an enumeration of all closed nowhere dense sets in $2^\kappa$, and $[\kappa]^\kappa=\{X_\alpha\colon \alpha<2^\kappa\}$. We construct a~sequence $\langle x_\alpha\rangle_{\alpha<2^\kappa}\in (2^\kappa)^{2^\kappa}$ by induction. For $\alpha<2^\kappa$, using Proposition~\ref{mmrn}, choose a~$\kappa$-meagre set $B_\alpha\subseteq 2^\kappa$ such that $F_\alpha+(2^\omega\setminus B_\alpha)$ is $\kappa$-Ramsey null. Pick any
\[x_\alpha\in \left\{\chi_S\colon S\in [X_\alpha]^\kappa\right\}\setminus \bigcup_{\beta<\alpha} \left(F_\beta+(2^\omega\setminus B_\beta)\right).\]
Such $x_\alpha$ exists, because $\cov(\kappa-CR_0)= 2^\kappa$. Indeed, if $X_\alpha=\{\eta_\xi\colon \xi<\kappa\}$, and $\left\{\chi_S\colon S\in [X_\alpha]^\kappa\right\}$ is covered by $\kappa-CR_0$ sets $\langle R_\beta\rangle_{\beta<\alpha}$, then $2^\kappa$ is covered by $\langle R'_\beta\rangle_{\beta<\alpha}$, where 
$R'_\beta=\{x'\colon x\in R_{\beta}\}$ with $x'(\xi)=x(\eta_\xi)$. Moreover, for every $\beta<\alpha$, $R'_\beta$ is $\kappa-CR_0$. 

Let $A=\{x_\alpha\colon \alpha<2^\kappa\}$. Obviously, $A$ is not $\kappa$-Ramsey null, because for all $S\in [\kappa]^\kappa$, there exists $S'\in [S]^\kappa$ such that $\chi_{S'}\in A$. 

Moreover, if $F$ is nowhere dense, then let $\alpha<2^\kappa$ be such that $F\subseteq F_\alpha$. For every $\beta>\alpha$, $x_\beta\notin F_\alpha+(2^\omega\setminus B_\alpha)$, thus by Lemma~\ref{lem-mm}, $x_\beta+F_\alpha\subseteq B_\alpha$. Hence,
\[A+F\subseteq A+F_\alpha=\bigcup_{\beta\leq\alpha} (x_\beta+F_\alpha)\cup \bigcup_{\alpha<\beta<2^\kappa} (x_\beta+F_\alpha)=\bigcup_{\beta\leq\alpha} (x_\beta+F_\alpha)\cup B_\alpha,\]
which is $\kappa$-meagre, since $\add(\MM_\kappa)=2^\kappa$. \hfill $\square$

We leave the following question open.

\begin{prob}
Does the hypothesis of Theorem~\ref{manrams} hold under weaker assumptions?
\end{prob}

\subsection{$\kappa$-$v_0$-Sets}

In \cite{mkantw:ssrtfn}, one can find a definition of a \df{$v_0$-set}.  A~set $A\subseteq 2^\omega$ has this property, if for every Silver perfect set $P$, there exists a~Silver perfect set $Q\subseteq P$ with $Q\cap A=\zbp$.

We propose a similar notion. A~perfect set $P\subseteq 2^\kappa$ is called a~\df{$\kappa$-Silver perfect} if for all $\alpha<\kappa$ and any $i\in\{0,1\}$, there exists $s\in 2^\alpha\cap T_P$ such that $s^\frown i\in T_P$ if and only if for all $s\in 2^\alpha\cap T_P$, $s^\frown i\in T_P$.

We define two classes of sets depending of whether we take into account all perfect sets or $\kappa$-perfect sets only. Namely, a~set $A\subseteq 2^\kappa$ is a~\df{$\kappa$-$v_0$-$\kappa$-set} if for all $\kappa$-Silver $\kappa$-perfect set $P\subseteq 2^\kappa$, there exists a~$\kappa$-Silver $\kappa$-perfect set $R\subseteq P$ such that $A\cap R=\zbp$. A~set $A\subseteq 2^\kappa$ is a~\df{$\kappa$-$v_0$-set} if for all $\kappa$-Silver perfect set $P\subseteq 2^\kappa$, there exists a~$\kappa$-Silver perfect set $R\subseteq P$ such that $A\cap R=\zbp$. The notion of $\kappa$-$v_0$-sets was considered in \cite{gl:gssm}. We study the relation between this notion and other notions of special subsets of $2^\kappa$. 

The following observation in the case of $\kappa=\omega$ was made in the proof of \cite{mkantw:ssrtfn}[Theorem~2.1].

\begin{prop}[$\omega$: \cite{mkantw:ssrtfn}]\label{comv0}
Assume that $\kappa$ is a~strongly inaccessible cardinal. Then every $\kappa$-comeagre subset of $2^\kappa$ contains a~$\kappa$-Silver $\kappa$-perfect set. 
\end{prop}

Proof: Let $A\subseteq 2^\kappa$ be $\kappa$-meagre. By Proposition~\ref{char-meagre} we get an~increasing sequence $\langle \xi_\alpha\rangle_{\alpha<\kappa}\in \kappa^\kappa$ and $z\in 2^\kappa$ such that 
$A\subseteq \left\{x\in 2^\kappa\colon \exists_{\beta<\kappa}\forall_{\beta<\gamma<\kappa} \exists_{\xi_{\gamma}\leq\xi<\xi_{\gamma+1}} z(\xi)\neq x(\xi)\right\}$.

Let 
$R=\{x\in 2^\kappa \colon \forall_{\alpha\in \Lim} \forall_{\xi_{\alpha}\leq\xi<\xi_{\alpha+1}} x(\xi)=z(\xi)\}$.
Then $R\subseteq 2^\kappa\setminus A$, and $R$ is a~$\kappa$-Silver perfect set. \hfill $\square$

\begin{cor}[$\omega$: \cite{mkantw:ssrtfn}, Theorem~2.1]\label{pkmkv0}
Assume that $\kappa$ is a~strongly inaccessible cardinal. Then every $\kappa$-perfectly $\kappa$-meagre set in $2^\kappa$ is a~$\kappa$-$v_0$-$\kappa$-set.
\end{cor}

Proof: Notice that for every $\kappa$-Silver $\kappa$-perfect set $P\subseteq 2^\kappa$, there exists a~natural homeomorphism $h\colon P\to 2^\kappa$ such that $R\subseteq 2^\kappa$ is a~$\kappa$-Silver $\kappa$-perfect set if an only if $h^{-1}[R]$ is $\kappa$-Silver $\kappa$-perfect. The corollary thus follows from Proposition~\ref{comv0}. 

\ \hfill $\square$

The following problem remains open.

\begin{prob}
Does the hypothesis of Corollary~\ref{pkmkv0} hold under weaker assumptions?
\end{prob}

On the other hand, regardless of properties of $\kappa$, we prove the following fact.

\begin{prop}[$\omega$: \cite{mkantw:ssrtfn}, Corollary~2.5]\label{snv0}
Every $\kappa$-strongly null set in $2^\kappa$ is a~$\kappa$-$v_0$-$\kappa$-set.
\end{prop}

Proof: Let $P\subseteq 2^\kappa$ be a~$\kappa$-Silver $\kappa$-perfect set, and $A\in\SN_\kappa$. Let 
$S=\{\len(s)\colon s\in \Spl(T_P)\}$.
Let $b\colon \kappa\times \{0,1\}\to S$ be a~bijection, and let $X=b[\kappa\times\{0\}]$. Let $\langle x_\alpha\rangle_{\alpha\in X}\in (2^\kappa)^X$ be such that $A\subseteq \bigcup_{\alpha\in X} [x_\alpha\obc \alpha+1]$. Then
$R=\left\{x\in P\colon \forall_{\alpha\in X} x(\alpha)=x_\alpha(\alpha)+1\right\}$
is a~$\kappa$-Silver $\kappa$-perfect set such that $R\subseteq P$, and $R\cap A=\zbp$. \hfill $\square$   
 
Analogously to Question~\ref{qrams}, we may ask the following question.

\begin{prob}
Are the ideals of $\kappa$-$v_0$-sets and $\kappa$-$v_0$-$\kappa$-sets  $\kappa^+$-complete?
\end{prob}

\subsection{$\kappa$-$l_0$-Sets}  

Recall that a~set $A\subseteq \omega^\omega$  is a~\df{$l_0$-set} if for every Laver perfect set $P$, there exists a~ Laver  perfect set $Q\subseteq P$ with $Q\cap A=\zbp$ (see \cite{mktw:ssrtfn}).

A perfect pruned tree $T\subseteq \kappa^{<\kappa}$ is called a~\df{$\kappa$-Laver perfect tree} if there exists $s\in T$ such that for all $t\in T$, either $t\subseteq s$, or 
$\left|\left\{\alpha<\kappa\colon t^\frown  \alpha\in T\right\}\right|=\kappa$.
A perfect set $P$ is \df{$\kappa$-Laver}, if $T_P$ is a $\kappa$-Laver perfect tree.

As before, we define two classes of sets depending of whether we take into account all perfect sets or $\kappa$-perfect sets only. A~set $A\subseteq \kappa^\kappa$ is \df{$\kappa$-$l_0$-$\kappa$-set}  if for every $\kappa$-Laver  $\kappa$-perfect set $P$, there exists a~$\kappa$-Laver  $\kappa$-perfect set $Q\subseteq P$ such that $Q\cap A=\zbp$.   
A~set $A\subseteq \kappa^\kappa$ is \df{$\kappa$-$l_0$-set}  if for every $\kappa$-Laver  perfect set $P$, there exists a~$\kappa$-Laver  perfect set $Q\subseteq P$ such that $Q\cap A=\zbp$.

\begin{theo}[$\omega$: \cite{mktw:ssrtfn}, Theorem~4.2]\label{snl0}
Every $\kappa$-strongly null set in $\kappa^\kappa$ is a~$\kappa$-$l_0$-$\kappa$-set.
\end{theo}

Proof: Let $T\subseteq \kappa^{<\kappa}$ be a $\kappa$-perfect $\kappa$-Laver tree, and suppose that $A\subseteq \kappa^\kappa$ is a $\kappa$-strongly null set. Let $t_0\in T$ be such that for every $s\in T$ with $t_0\subseteq s$, $\left|\left\{\alpha<\kappa\colon t_0\,^\frown  \alpha\in T\right\}\right|=\kappa$.

 Let $I=\{\alpha<\kappa\colon \len(t)<\alpha\}\setminus \Lim$, and assume that $\langle s_{\alpha}\rangle_{\alpha\in I}$ satisfies $A=\bigcup_{\alpha\in I} [s_\alpha]$,
 where for all $\alpha\in I$, $s_\alpha\in 2^\alpha$.
 
We construct tree $T'\subseteq \kappa^{<\kappa}$ in the following way. Let $T'\cap \kappa^{\leq\len(t_0)}=T\cap \kappa^{\leq\len(t_0)}$, and assume that $\alpha<\kappa$ is such that $\alpha> \len(t_0)$, and $t\in T\cap \kappa^\alpha$. Then put 
\[T'\cap \{s\in \kappa^{\alpha+1}\colon t\subseteq s\}= T \cap \{s\in \kappa^{\alpha+1}\colon t\subseteq s\}\setminus \{s_\alpha\}.\]
For limit $\beta<\kappa$ with $\beta>\len(t_0)$, let $t\in T'$ if and only if $t\obc\alpha\in T'$ for every $\alpha<\beta$. 

Since $\kappa$ is regular, $T'$ is a $\kappa$-perfect $\kappa$-Laver tree, and $[T']\subseteq [T]\setminus A$. \hfill $\square$

The following problem remains open.

\begin{prob}
Does the hypothesis of Theorem~\ref{snl0} hold under weaker assumptions?
\end{prob}

Analogously to Question~\ref{qrams}, we may also ask the following question.

\begin{prob}
Are the ideals of $\kappa$-$l_0$-sets and $\kappa$-$l_0$-$\kappa$-sets  $\kappa^+$-complete?
\end{prob}

\subsection{$\kappa$-$m_0$-Sets} \label{s-m0}

In \cite{mktw:ssrtfn}, the authors consider the notion of \df{$m_0$-set}. A~set $A\subseteq \omega^\omega$ is an $m_0$-set if for every Miller perfect set $P$, there exists a~ Miller  perfect set $Q\subseteq P$ with $Q\cap A=\zbp$.

We define an analogous notion. A perfect pruned tree $T\subseteq \kappa^{<\kappa}$ is called a~\df{$\kappa$-Miller perfect tree} if for every $s\in T$ there exists $t\in T$ such that  $s\subseteq t$, and 
$\left|\left\{\alpha<\kappa\colon t^\frown \alpha\in T\right\}\right|=\kappa$.
A perfect set $P$ is \df{$\kappa$-Miller}, if $T_P$ is a $\kappa$-Miller perfect tree. 

Similarly as before, we define two classes of sets depending of whether we take into account all perfect sets or $\kappa$-perfect sets only, a~set $A\subseteq 2^\kappa$ is \df{$\kappa$-$m_0$-$\kappa$-set} if for every $\kappa$-Miller $\kappa$-perfect set $P$, there exists a~ $\kappa$-Miller $\kappa$-perfect set $Q\subseteq P$ such that $Q\cap A=\zbp$.  
A~set $A\subseteq 2^\kappa$ is \df{$\kappa$-$m_0$-set} if for every $\kappa$-Miller perfect set $P$, there exists a~ $\kappa$-Miller perfect set $Q\subseteq P$ such that $Q\cap A=\zbp$.

\begin{theo}[$\omega$: \cite{mktw:ssrtfn}, Theorem~3.2]
Every $\kappa$-perfectly $\kappa$-meagre set in $\kappa^\kappa$ is a~$\kappa$-$m_0$-$\kappa$-set.
\end{theo}

Proof: Let $P$ be a $\kappa$-perfect $\kappa$-Miller set, and $A\in\PkM_\kappa$. Let $h\colon P\to \kappa^\kappa$ be the standard homeomorphism induced by the order isomorphism between $P$ and $\kappa^\kappa$. Notice also that under this homeomorphism, if $Q\subseteq \kappa^\kappa$ is a $\kappa$-perfect $\kappa$-Miller set, then $h^{-1}[Q]$ is a $\kappa$-Miller $\kappa$-perfect set as well. 

Thus, $B=h[A\cap P]$ is $\kappa$-meagre. Let $B=\bigcup_{\alpha<\kappa} G_\alpha$ with $G_\alpha$ nowhere dense closed for every $\alpha<\kappa$. 

We choose by induction $\left<t_s\right>_{s\in \kappa^{<\kappa}}$ such that $t_s \in \kappa^{<\kappa}$, for every $s\in \kappa^{<\kappa}$, and
\[\left|\{\alpha<\kappa \colon \exists_{s'\in \kappa^{\kappa}}\colon t_s\,^\frown \alpha\subseteq t_{s'}\}\right|=\kappa,\]
 and for $s,s'\in \kappa^{<\kappa}$, $s\subsetneq s'$ if and only if $t_s\subsetneq t_{s'}$. Indeed, let $t_\zbp$ be such that $[t_\zbp]\cap G_0=\zbp$
Then, given $t_s$, $s\in \kappa^\alpha$, let $t'_s\supsetneq t_s$ be such that $[t'_s]\cap G_{\alpha+1}= \zbp$. Set $t_{s^\frown \xi}=t_s^{\prime\frown}\xi$, for all $\xi<\kappa$. For limit $\beta<\kappa$, and $s\in \kappa^\beta$, put $t'_s=\bigcup_{\alpha<\beta}t_{s\obc\alpha}$. Let $t_s\supsetneq t'_s$ be such that $[t_s]\cap G_\beta=\zbp$. Finally, define $T=\bigcup_{\alpha<\kappa}\{t\in \kappa^{<\kappa}\colon t\subseteq t_{s}, s\in \kappa^{\alpha}\}$. Obviously, $T$ is a~$\kappa$-perfect $\kappa$-Miller tree, so $P'=[T]_{\kappa}$ is a~$\kappa$-perfect  $\kappa$-Miller subset of $\kappa^\kappa\setminus B$.

Thus, there exists a $\kappa$-perfect $\kappa$-Miller $Q\subseteq P\setminus A$. \ \hfill $\square$

Next, we will need the following lemma.

\begin{lemma}\label{lem-fod}
Let $\oF\subseteq \kappa^\kappa$ with $|\oF|\leq \kappa$. Then there exists a function $g\colon \kappa\to \kappa$ such that for every injection $i\colon \kappa\to \kappa$ and every $f\in \oF$, $g\circ i\not\leq^\kappa f$.
\end{lemma}

Proof: Obviously, by Fodor's Lemma there exists cofinally many $\alpha<\kappa$ such that $i(\alpha)\geq \alpha$. Let $g$ be an increasing function such that for every $f\in \oF$, $g\geq^\kappa f$. But then if for an injection $i$, $f\geq^\kappa g\circ i$, there exists cofinally many $\alpha<\kappa$ such that
$f(\alpha)\geq g(i(\alpha))\geq g(\alpha),$
which is impossible. \hfill $\square$

\begin{theo}[$\omega$: \cite{mktw:ssrtfn}, Theorem~4.3]\label{m0th}
Every $\kappa$-strongly null set in $\kappa^\kappa$ is a~$\kappa$-$m_0$-set.
\end{theo}

Proof: Let $T\subseteq \kappa^\kappa$ be a $\kappa$-Miller perfect tree. 
Let 
$\Spl_\kappa(T)=\left\{s\in T\colon  \left|\left\{\alpha\in \kappa \colon s^\frown \alpha\in T\right\}\right|=\kappa\right\}$.
For $t\in T$, let $s(t)\in \Spl_\kappa(T)$ be any minimal element of 
$\left\{s\in T\colon t\subseteq s \land s\in \Spl_\kappa(T)\right\}$, and let 
$S(t)=\{s(s(t)^\frown \alpha)\colon \alpha<\kappa \land s(t)^\frown \alpha \in T\}$.
Also let $\{s(t,\alpha)\colon \alpha<\kappa\}$ be an enumeration of $S(t)$, and for $\beta<\kappa$ let $f_{t,\beta}\in \kappa^\kappa$ with $f_{t,\beta}(\alpha)=\len(s(t,\beta+\alpha))$.

Moreover, for $\xi\in \kappa\cap \Lim$ let 
\[L_{\xi}(t)=\left\{s\in T\colon \exists_{\langle s_\beta\rangle_{\beta<\xi}\in (\Spl_\kappa(T))^\xi} \left(\forall_{\beta<\xi} t\subseteq s_\beta \land s\obc \len(s_\beta)=s_\beta  \right) \land \left(\forall_{\beta<\gamma<\xi} s_\beta \subsetneq s_\gamma \right)\right\}.\]
Put $S_\xi(t)=\{s(u)\colon u\in L_\xi(t)\}$. Finally, let $\{s(t,\xi,\alpha)\colon \alpha<\kappa\}$ be an enumeration of $S_\xi(t)$, and for $\beta<\kappa$ let $f_{t,\xi,\beta}\in \kappa^\kappa$ with $f_{t,\xi,\beta}(\alpha)=\len(s(t,\xi,\beta+\alpha))$.

Define $\oF= \{f_{t,\beta}\colon t\in T, \beta<\kappa\}\cup \{f_{t,\xi,\beta}\colon t\in T, \xi,\beta<\kappa\}$. By Lemma~\ref{lem-fod}, there exists $g\in \kappa\to\kappa$ such that for all $t\in T$,  $\xi,\beta<\kappa$, and every injection $i\colon \kappa\to \kappa$, $g\circ i\not\leq^\kappa f_{t,\beta}$ and $g\circ i\not\leq^\kappa f_{t,\xi, \beta}$. Notice that we may assume that $g(0)>\len(s(\zbp))$. 

Let $A$ be a $\kappa$-strongly null set, and let $\langle s_\alpha\rangle_{\alpha<\kappa}\in (\kappa^{<\kappa})^\kappa$ be such that $s_\alpha \in 2^{g(\alpha)}$, for all $\alpha<\kappa$, and
$A\subseteq \bigcup_{\alpha<\kappa}[s_\alpha]$.

Notice that for every $t\in T$ such that there is no $\alpha<\kappa$ with $s_\alpha\subseteq s(t)$,  
$|S(t)\setminus \{s\in T\colon \exists_{\alpha<\kappa} s_\alpha\subseteq s\}|=\kappa$.

Assume otherwise that there exists $t\in T$ such that there is no $\alpha<\kappa$ with $s_\alpha\subseteq s(t)$,  but there exists $\beta_0<\kappa$ such that for all $\beta_0\leq\alpha<\kappa$ there exists $\xi_\alpha<\kappa$ such that $s_{\xi_{\alpha}}\subseteq s_{t,\alpha}$. Let $i\colon \kappa\to \kappa$ be defined as $i(\alpha)=\xi_{\beta_0+\alpha}$. Since there is no $\alpha<\kappa$ with $s_\alpha\subseteq s(t)$, $i$ is an injection. Thus, we obtain for all $\alpha<\kappa$, 
\[g(i(\alpha))=g(\xi_{\beta_0+\alpha})=\len(s_{\xi_{\beta_0+\alpha}})\leq \len(s(t,\beta_0+\alpha))=f_{t,\beta_0}(\alpha),\]
which is a contradiction.

A very similar argument leads to the conclusion that if $t\in T$ is such that there is no $\alpha<\kappa$ with $s_\alpha\subseteq s(t)$, and  $\xi\in \kappa\cap \Lim$, then  
$|S_\xi(t)\setminus \{s\in T\colon \exists_{\alpha<\kappa} s_\alpha\subseteq s\}|=\kappa$.

Now, consider the set
 $T'=\{t\in T\colon \forall_{\alpha<\kappa} s_\alpha\not\subseteq s(t)\}$.
It follows that $T'$ is a pruned perfect tree, which is $\kappa$-Miller. Thus $[T']$ is a $\kappa$-Miller perfect set, and $[T']\cap A=\zbp$. \hfill $\square$

We leave the following problem open.

\begin{prob}
Is every $\kappa$-strongly null set in $\kappa^\kappa$ a~$\kappa$-$m_0$-$\kappa$-set.
\end{prob}

Similarly to Question~\ref{qrams}, we may ask the following question.

\begin{prob}
Are the ideals of $\kappa$-$m_0$-sets and $\kappa$-$m_0$-$\kappa$-sets  $\kappa^+$-complete?
\end{prob}

\providecommand{\WileyBibTextsc}{}
\let\textsc\WileyBibTextsc
\providecommand{\othercit}{}
\providecommand{\jr}[1]{#1}
\providecommand{\etal}{~et~al.}


\end{document}